\documentclass{article}
\usepackage{amsmath,amssymb,amsthm}
\begin{document}
\newtheorem*{thm*}{Theorem}
\newtheorem*{proposition*}{Proposition}
\newtheorem*{corollary*}{Corollary}



\def\al{\alpha}
\def\be{\beta}
\def\ga{\gamma}
\def\Ga{\Gamma}
\def\la{\lambda}
\def\om{\omega}
\def\Om{\Omega}

\def\C{{\mathbb C}}
\def\R{{\mathbb R}}
\def\Z{{\mathbb Z}}
\def\Y{{\mathbb Y}}

\def\one{{\mathbf1}}
\def\Pz{\widetilde{\mathcal P}_z}
\def\Cal#1{{\mathcal{#1}}}
\def\n{{(n)}}

\def\sgn{\operatorname{sgn}}
\def\diag{\operatorname{diag}}
\def\Ex{\operatorname{Ex}}
\def\bS{\mathfrak S}
\def\Prob{{\operatorname{Prob}}}
\def\Conf{{\operatorname{Conf}}}
\def\Tr{{\operatorname{Tr}}}
\def\ctg{\operatorname{ctg}}


\title{An
introduction to harmonic analysis on the infinite symmetric group}
\author{Grigori Olshanski}




\maketitle

\section*{Introduction}

The aim of the present survey paper is to provide an accessible
introduction to a new chapter of representation theory --- harmonic
analysis for {\it noncommutative groups with infinite--dimensional
dual space.}

I omitted detailed proofs but tried to explain the main ideas of the theory and
its connections with other fields. The fact that irreducible representations of
the groups in question depend on infinitely many parameters leads to a number
of new effects which never occurred in conventional noncommutative harmonic
analysis. A link with stochastic point processes is especially emphasized.

The exposition focuses on a single group, the infinite symmetric group
$S(\infty)$. The reason is that presently this particular example is worked out
the most. Furthermore, $S(\infty)$ can serve as a very good model for more
complicated groups like the infinite--dimensional unitary group $U(\infty)$.

The paper is organized as follows. In \S1, I explain what is the problem of
harmonic analysis for $S(\infty)$. \S\S2--5 contain the necessary preparatory
material. In \S6, the main result is stated. It was obtained in a cycle of
papers by Alexei Borodin and myself. In \S7, the scheme of the proof is
outlined. The final \S8 contains additional comments and detailed references.

This paper is an expanded version of lectures I gave at the Euler Institute,
St.--Petersburg, during the NATO ASI Program ``Asymptotic combinatorics with
applications to mathematical physics''. I also partly used the material of my
lectures at the Weizmann Institute of Science, Rehovot. I am grateful to
Anatoly Vershik, Amitai Regev, and Anthony Joseph for warm hospitality in
St.~Petersburg and Rehovot, and to Vladimir Berkovich for taking notes of my
lectures at the Weizmann. Finally, I would like to thank Alexei Borodin for
cooperation and help.

\section{Virtual permutations and generalized regular representations}

\subsection{The Peter--Weyl theorem}

Let $\Cal K$ be a compact group, $\mu$ be the normalized Haar measure on
$\Cal K$ (i.e., $\mu(\Cal K)=1$), and $H$ be the Hilbert space
$L^2(K,\mu)$. The
group $\Cal G=\Cal K\times \Cal K$ acts on $\Cal K$ on the right as
follows: if $g=(g_1,g_2)\in G$ and $x\in \Cal K$, then $x\cdot
g=g_2^{-1}x g_1$. This action gives rise to a unitary representation
$T$ of $\Cal G$ on $H$:
$$
(T(g)f)(x)=f(x\cdot g), \qquad f\in H, \quad g\in G.
$$
It is called the {\it biregular representation.\/} Let $\widehat{\Cal K}$ be
the set of equivalence classes of irreducible representations of $\Cal K$.
Recall that all of them are finite dimensional and unitarizable. For
$\pi\in\widehat{\Cal K}$, let $\overline\pi$ denote the dual representation.
Since $\pi$ is unitary, $\overline\pi$ is obtained from $\pi$ by the
conjugation automorphism of the base field $\C$.

\newtheorem*{pwt}{Peter--Weyl's Theorem}
\begin{pwt}The biregular representation $T$
is equivalent to the direct sum of the irreducible representations of
$G$ of the form $\pi\otimes\overline\pi$,
$$
T\sim\bigoplus\limits_{\pi\in\widehat{\Cal K}}(\pi\otimes\overline\pi).
$$
\end{pwt}

This is one of the first results of {\it noncommutative harmonic
analysis.\/} The aim of noncommutative harmonic analysis can be
stated as {\it decomposing natural representations into irreducible
ones.\/} The biregular representation can be called a natural
representation because it is fabricated from the group itself in a
very natural way. The Peter--Weyl theorem serves as a guiding example
for more involved theories of noncommutative harmonic analysis.

\subsection{The infinite symmetric group}

Let $S(n)$ be the symmetric group of degree $n$, i.e., the group of
permutations of the set $\{1,\dots,n\}$. By the very definition,
$S(n)$ acts on $\{1,\dots,n\}$. The stabilizer of $n$ is canonically
isomorphic to $S(n-1)$, which makes it possible to define, for any
$n=2,3,\dots$, an embedding $S(n-1)\to S(n)$. Let $S(\infty)$ be the
inductive limit of the groups $S(n)$ taken with respect to these
embeddings. We call $S(\infty)$ the {\it infinite symmetric group.}

Clearly, $S(\infty)$ is a countable, locally finite group. It can be
realized as the group of all {\it finite\/} permutations of the set
$\{1,2,\dots\}$\,.

\subsection{The biregular representation for $S(\infty)$}

The definition of a biregular representation given in \S1.1
evidently makes sense for the group $S(\infty)$. Namely, set
$\Cal K=S(\infty)$, $\Cal G=S(\infty)\times S(\infty)$, and take as
$\mu$ the counting measure on $S(\infty)$. Then the
unitary representation $T$ of the group $S(\infty)\times S(\infty)$
in the Hilbert space $L^2(S(\infty),\mu)$ is defined by exactly the
same formula as in \S1.1.

\begin{proposition*}
The biregular representation $T$ of
the group $S(\infty)\times S(\infty)$ is
irreducible.
\end{proposition*}

{\renewcommand{\proofname}{Sketch of proof}
\begin{proof}
Let $\diag(S(\infty))$ be the image of
$S(\infty)$ under the diagonal embedding $S(\infty)\to
S(\infty)\times S(\infty)$. The Dirac function
$\delta_e$ is a unique (up to a constant factor)
$\diag(S(\infty))$--invariant vector in the space of $T$. This
follows from the fact that all
conjugacy classes in $S(\infty)$, except for $\{e\}$, are infinite. On
the other hand, $\delta_e$ is a cyclic vector, i.e., it generates under
the action of $S(\infty)\times S(\infty)$ a dense subspace in
$L^2(S(\infty)$). It follows that there is no proper closed
$S(\infty)\times S(\infty)$--invariant subspace.
\end{proof}}

Thus, in the case of the group $S(\infty)$, the naive analog of the biregular
representation is of no interest for harmonic analysis. We will explain how to
modify the construction in order to get interesting representations.

{}From now on we are using the notation
$$
G=S(\infty)\times S(\infty), \qquad K=\diag(S(\infty)).
$$
We call $G$ the {\it infinite bisymmetric group.\/}

\subsection{Virtual permutations}

Note that in the construction of \S1.1, the group $\Cal K$
plays two different roles: it is the carrier of a Hilbert space of
functions and it acts (by left and right shifts) in this space. The idea
is to separate these two roles. As the carrier of
a Hilbert space we will use a remarkable compactification $\bS$ of
$S(\infty)$. It is not a group but still a $G$--space, which is sufficient
for our purposes.

For any $n\ge2$, we define a projection $p_n: S(n)\to S(n-1)$ as
removing the element $n$ from the cycle containing it. That is, given
a permutation $\sigma\in S(n)$, if $n$ is fixed under $\sigma$ then
$p_n(\sigma)=\sigma$, and if $n$ enters a nontrivial cycle
$(\cdots\to i\to n\to j\to\cdots)$ then we simply replace this cycle by
$(\cdots\to i\to j\to\cdots)$. We call $p_n$ the {\it canonical
projection}.

\begin{proposition*} The canonical projection $p_n:S(n)\to S(n-1)$
commutes with the left and right shifts by the elements of $S(n-1)$.
Moreover, for $n\ge5$ it is the only map $S(n)\to S(n-1)$ with such a
property.
\end{proposition*}

Let $\bS$ be the projective limit of the finite sets $S(n)$ taken
with respect to the canonical projections. Any point $x\in\bS$ is a
collection $(x_n)_{n\ge1}$ such that $x_n\in S(n)$ and
$p_n(x_n)=x_{n-1}$. For any $m$, we identify $S(m)$ with the subset
of those points $x=(x_n)$ for which $x_n\in S(m)$ for all $n\ge m$.
This allows us to embed $S(\infty)$ into $\bS$.

We equip $\bS$ with the projective limit topology. In this way we get
a totally disconnected compact topological space. We call it the {\it
space of virtual permutations.}

The image of $S(\infty)$ is dense in $\bS$. Hence, $\bS$ is a
{\it compactification\/} of the discrete space $S(\infty)$.

There exists an action of the group $G$ on the space $\bS$ by
homeomorphisms extending the action of $G$ on $S(\infty)$. Such an
action is unique.

There are several different realizations of the space $\bS$. One of them looks
as follows. Set $I_n=\{0,\dots,n-1\}$. There exists a bijection
$$
\bS\to I:=I_1\times I_2\times\dots, \qquad
x=(x_1,x_2,\dots)\mapsto(i_1,i_2,\dots)
$$
such that $i_n=0$, if $x_n(n)=n$, and $i_n=j$, if $x_n(n)=j<n$.
This bijection is a homeomorphism (here we equip the product space
$I$ with the product topology). It gives rise, for every $n\ge1$, to
a bijection $S(n)\to I_1\times\dots\times I_n$. In this realization,
the canonical projection $p_n:S(n)\to S(n-1)$ turns into the natural
projection $I_1\times\dots\times I_n\to I_1\times\dots\times I_{n-1}$.

\subsection{Ewens' measures on $\bS$}

Let $\mu_1^{(n)}$ be the normalized Haar measure on $S(n)$. Its
pushforward under the canonical projection $p_n$ coincides with the
measure $\mu_1^{(n-1)}$, because $p_n$ commutes with the left (and
right) shifts by elements of $S(n-1)$. Thus, the measures
$\mu_1^{(n)}$ are pairwise {\it consistent\/} with respect to the
canonical projections. Hence, we can define their projective limit,
$\mu_1=\varprojlim\mu_1^{(n)}$, which is a probability
measure on $\bS$.

The measure $\mu_1$ is invariant under the action of $G$, and it is
the only probability measure on $\bS$ with this property. Thus,
viewing $\bS$ as a substitute of the group space, we may view $\mu_1$
as a substitute of the normalized Haar measure.

Now we define a one--parameter family of probability measures
containing the measure $\mu_1$ as a particular case.

For $t\ge0$, let $\mu^{(n)}_t$ be the following measure on $S(n)$:
$$
\mu^{(n)}_t(x)=\frac{{t^{[x]-1}}{(t+1)(t+2)\cdot\dots\cdot(t+n-1)}}\,,
$$
where $[x]=[x]_n$ is the number of cycles of $x$ in $S(n)$. If
$t=1$ then this reduces to above definition of the measure $\mu_1^{(n)}$.

\begin{proposition*} {\rm(i)} $\mu^{(n)}_t$ is a probability
measure on $S(n)$, i.e.,
$$
\sum_{x\in S(n)}t^{[x]}=t(t+1)\cdot\dots\cdot(t+n-1).
$$

{\rm(ii)} The measures $\mu^{(n)}_t$ are pairwise consistent with
respect to the canonical projections.

{\rm(iii)} The pushforward of $\mu_t^{(n)}$ under the bijective map
$S(n)\to I_1\times\dots\times I_n$ of \S1.4 is the product measure
$\nu^{(1)}_t\times \dots\times\nu^{(n)}_t$, where, for any $m$,
$\nu_t^{(m)}$ is the following probability measure on $I_m$:
$$
\nu^{(m)}_t(i)=\begin{cases} \frac t{t+m-1}, & i=0\\
\frac1{t+m-1}, & i=1,\dots,m-1. \end{cases}
$$
\end{proposition*}

\begin{proof} (i) Induction on $n$.
Assume that the equality in question holds for $n-1$. Notice that
$[p_n(x)]_{n-1}$ is equal to $[x]_n$ when $x\not\in S(n-1)$, and to
$[x]_n-1$ when $x\in S(n-1)$. We have
\begin{align*}
\sum_{x\in S(n)}t^{[x]}&=\sum_{y\in S(n-1)}
\sum_{p_n(x)=y}t^{[x]}=\sum_{y\in S(n-1)}\left(t\cdot t^{[y]} +
(n-1)t^{[y]}\right)\\
&=t\cdot t(t+1)\cdot\dots\cdot(t+n-2) +
(n-1)t(t+1)\cdot\dots\cdot(t+n-2)\\
&=t(t+1)\cdot\dots\cdot(t+n-1).
\end{align*}

(ii) We have to verify that for every $y\in S(n-1)$
$$
\frac{{t^{[y]-1}}{(t+1)\cdot\dots\cdot(t+n-2)}}= \sum_{p_n(x)=y}
\frac{{t^{[x]-1}}{(t+1)\cdot\dots\cdot(t+n-1)}}\ .
$$
It is precisely what is done in the proof of (i).

(iii) This follows from the fact that, under the bijection
$x\mapsto(i_1,\dots,i_n)$ between $S(n)$ and $I_1\times\dots\times
I_n$, the number of zeros in $(i_1,\dots,i_n)$ equals $[x]$.
\end{proof}

The consistency property makes it possible to define, for any
$t\ge0$, a probability measure $\mu_t=\varprojlim\mu^{(n)}_t$ on
$\bS$. This measure is invariant under the diagonal subgroup $K$
but is not $G$--invariant (except the case $t=1$).  As $t\to\infty$,
$\mu_t$ tends to the Dirac measure at $e\in S(\infty)\subset\bS$. Let
us denote this limit measure by $\mu_\infty$.

Following S.~V.~Kerov, we call the measures $\mu_t$ the {\it Ewens
measures.\/} The next claim gives a characterization of the family
$\{\mu_t\}$.

\begin{proposition*} The measures $\mu_t$, where $0\le t\le\infty$,
are exactly those probability measures on $\bS$ that are
$K$--invariant and correspond to product measures on
$I_1\times I_2\times\dots$\,.
\end{proposition*}

\subsection{Transformation properties of the Ewens measures}

Recall that $[\sigma]_n$ denotes the number of cycles of a permutation
$\sigma\in S(n)$.

\begin{proposition*} {\rm(i)} For any $x=(x_n)\in\bS$ and $g\in G$, the
quantity $[x_n\cdot g]_n-[x_n]_n$ does not depend on $n$ provided
that $n$ is large enough.

{\rm(ii)} Denote by $c(x,g)$ the stable value of this quantity. The
function $c(x,g)$ is an additive cocycle with values in $\Z$, that is,
$$
c(x,gh)=c(x\cdot g,h)+c(x,g), \qquad x\in\bS, \quad g,h\in G.
$$
\end{proposition*}

Recall that a measure is called {\it quasi--invariant\/} under a
group of transformations if, under the shift by an arbitrary element
of the group, the measure is transformed to an equivalent measure.

\begin{proposition*} Assume $t\in(0,+\infty)$.

{\rm(i)} The measure $\mu_t$ is quasi--invariant under the action of
the group $G$.

{\rm(ii)} We have
$$
\frac{\mu_t(d(x\cdot g))}{\mu_t(dx)}=t^{c(x,g)},
$$
where the left--hand side is the Radon--Nikodym derivative.
\end{proposition*}

Note that $c(x,g)=1$ whenever $g\in K$. This agrees with the
fact that the measures are $K$--invariant.

\newsavebox{\olshsb}\sbox{\olshsb}{$T_z$}
\subsection{The representations \usebox{\olshsb}}

We start with a general construction of unitary representations
related to group actions on measure spaces with cocycles.

Assume we are given a space $\Cal S$ equipped with a Borel
structure (i.e., a distinguished sigma--algebra of sets), a discrete group
$\Cal G$ acting on $\Cal S$ on the right and preserving the Borel
structure, and a Borel measure $\mu$, which is quasi--invariant under
$\Cal G$. A complex valued function $\tau(x,g)$ on $\Cal S\times\Cal
G$ is called a {\it multiplicative cocycle\/} if
$$
\tau(x,gh)=\tau(x\cdot g,h)\tau(x,g), \qquad x\in\Cal S,
\quad g,h\in\Cal G.
$$
Next, assume we are given a multiplicative cocycle $\tau(x,g)$ which is
a Borel function in $x$ and which satisfies the relation
$$
|\tau(x,g)|^2=\frac{\mu(d(x\cdot g))}{\mu(dx)}\,.
$$
Then these data allow us to define a unitary representation
$T=T_\tau$ of the group $\Cal G$ acting in the Hilbert space
$L^2(\Cal S,\mu)$ according to the formula
$$
(T(g)f)(x)=\tau(x,g)f(x\cdot g), \qquad f\in L^2(\Cal S,\mu), \quad x\in\Cal S,
\quad g\in\Cal G.
$$

Let $z\in\C$ be a nonzero complex number. We apply this general
construction for the space $\Cal S=\bS$, the group $\Cal G=G$, the
measure $\mu=\mu_t$ (where $t=|z|^2$), and the cocycle
$\tau(x,g)=z^{c(x,g)}$. All the assumptions above are satisfied, so
that we get a unitary representation $T=T_z$ of the group $G$.

Using a continuity argument it is possible to extend the definition
of the representations $T_z$ to the limit values $z=0$ and $z=\infty$
of the parameter $z$. It turns out that the representation $T_\infty$
is equivalent to the biregular representation of \S1.3. Thus, the
family $\{T_z\}$ can be viewed as a deformation of the biregular
representation.

We call the $T_z$'s the {\it generalized regular representations.\/}
These representations are reducible (with the only exception of
$T_\infty$). Now we can state the main problem that we address in
this paper.

\newtheorem*{proha}{Problem of harmonic analysis on $S(\infty)$}
\begin{proha}
Describe the decomposition of the generalized regular representations
$T_z$ into irreducibles ones.
\end{proha}

\section{Spherical representations and characters}

\subsection{Spherical representations}

By a {\it spherical representation\/} of the pair $(G,K)$ we mean a pair
$(T,\xi)$, where $T$ is a unitary representation of $G$ and $\xi$ is a
unit vector in the Hilbert space $H(T)$ such that:

(i) $\xi$ is $K$--invariant and

(ii) $\xi$ is cyclic, i.e., the span of the vectors of the
form $T(g)\xi$, where $g\in G$, is dense in $H(T)$.

We call $\xi$ the {\it spherical vector.\/} We call two spherical
representations $(T_1,\xi_1)$ and $(T_2,\xi_2)$ {\it equivalent\/} if
there exists an isometric isomorphism between their Hilbert spaces
which commutes with the action of $G$ and preserves the spherical
vectors. Such an isomorphism is unique within multiplication by a
scalar. The equivalence $(T_1,\xi_1)\sim(T_2,\xi_2)$ implies the
equivalence $T_1\sim T_2$ but the converse is not true in general.

The matrix coefficient $(T(g)\xi,\xi)$, where $g\in G$, is called the
{\it spherical function.\/} Two spherical representations are
equivalent if and only if their spherical functions coincide.

We aim to give an independent characterization of spherical functions
for $(G,K)$.

\subsection{Positive definite functions}

Recall that a complex--valued function $f$ on a group $\Cal G$ is
called {\it positive definite\/} if:

(i) $f(g^{-1})=\overline{f(g)}$ for any $g\in\Cal G$ and

(ii) for any finite collection $g_1,\dots,g_n$ of elements of $\Cal
G$, the $n\times n$ Hermitian matrix $[f(g_j^{-1}g_i)]$ is
nonnegative.

Positive definite functions on $\Cal G$ are exactly diagonal matrix
coefficients of unitary representations of $\Cal G$.

Now return to our pair $(G,K)$. The spherical functions for $(G,K)$
can be characterized as the {\it positive definite, $K$--biinvariant
functions on $G$, normalized at $e\in G$.}

\subsection{Characters}

Recall that the character of an irreducible representation $\pi$
of a compact group $\Cal K$ is the function $g\mapsto\chi^\pi(g)=
\Tr(\pi(g))$. If $\Cal K$ is noncompact, an irreducible representation
$\pi$ of $\Cal K$ is not necessarily finite dimensional, and so the
function $g\mapsto\Tr(\pi(g))$ does not make sense in general.
But it turns out that in certain cases the ratio
$$
\widetilde\chi^\pi(g)= \frac{\chi^\pi(g)}{\chi^\pi(e)}
$$
does make sense.

Let $\Cal K$ be an arbitrary group. A function on $\Cal K$ is said to be {\it
central\/} if it is constant on conjugacy classes. Denote by $\Cal X(K)$ the
set of central, positive definite, normalized functions on $\Cal K$ (if $\Cal
K$ is a topological group then we additionally require the functions to be
continuous). If $\varphi,\psi\in \Cal X(\Cal K)$, then for every $t\in[0,1]$
the function $(1-t)\varphi+t\psi$ is also an element of $\Cal X(\Cal K)$, i.e.,
$\Cal X(\Cal K)$ is a convex set.

Recall that a point of a convex set is called {\it extreme\/} if it
is not contained in the interior of an interval entirely contained in
the set. Let $\Ex(\Cal X(\Cal K))$ denote the subset of extreme
points of $\Cal X(\Cal K)$.

If the group $\Cal K$ is compact then the functions from $\Ex(\Cal
X(\Cal K))$ are exactly the {\it normalized\/} irreducible
characters $\widetilde\chi^\pi(g)$, where $\pi\in\widehat{\Cal K}$.
As for general elements of $\Cal X(\Cal K)$, they are (possibly
infinite) convex linear combinations of these functions.

In particular, if $\Cal K$ is finite then $\Cal X(\Cal K)$ is a
finite--dimensional simplex.

We will call the elements of $\Cal X(\Cal K)$ the {\it characters\/}
of $\Cal K$. The elements of $\Ex(\Cal X(\Cal K))$ will be called the
{\it extreme characters.\/} Notice that this terminology does not agree
with the conventional terminology of representation theory. However,
in the case of the group $S(\infty)$ this will not lead to a confusion.

\subsection{Correspondence between spherical representations of
$(G,K)$ and characters of $S(\infty)$}

There is a natural 1--1 correspondence between spherical functions
for $(G,K)$ and characters of $S(\infty)$. Specifically, given a
function $f$ on the group $G=S(\infty)\times S(\infty)$, let $\chi$
be the function on $S(\infty)$ obtained by restricting $f$ to the
first copy of $S(\infty)$. Then $f\mapsto\chi$ establishes a 1--1
correspondence between $K$--biinvariant functions on $G$ and central
functions on $S(\infty)$. Moreover, this correspondence preserves the
positive definiteness property. This implies that the equivalence
classes of spherical representations of $(G,K)$ are parametrized by
the characters of $S(\infty)$.

\begin{proposition*} Let $T$ be a unitary representation of $G$
and $H(T)^K$ be the subspace of $K$--invariant vectors in the Hilbert
space $H(T)$ of $T$.

If $T$ is irreducible then $H(T)^K$ has dimension 0 or 1. Conversely,
if the subspace $H(T)^K$ has dimension 1 then its cyclic span is an
irreducible subrepresentation of $T$.
\end{proposition*}

\begin{corollary*} For an irreducible spherical representation of
$(G,K)$, the spherical vector $\xi$ is defined uniquely, within a
scalar multiple, which does not affect the spherical function.
\end{corollary*}

A spherical function corresponds to an irreducible representation if
and only if the corresponding character is extreme. Thus, the
(equivalence classes of) irreducible spherical representations of
$(G,K)$ are parametrized by extreme characters of $S(\infty)$.

\subsection{Spectral decomposition}

\begin{proposition*} {\rm(i)} For any character
$\psi\in\Cal X(S(\infty))$, there exists a probability measure $P$ on
the set $\Ex(\Cal X(S(\infty)))$ of extreme characters such that
$$
\psi(\sigma)=\int_{\chi\in\Ex(\Cal X(S(\infty)))}
\chi(\sigma)P(d\chi),
\qquad \sigma\in S(\infty).
$$

{\rm(ii)} Such a measure is unique.

{\rm(iii)} Conversely, for any probability measure $P$ on the set of
extreme characters, the function $\psi$ defined by the above formula
is a character of $S(\infty)$.
\end{proposition*}

We call this integral representation the {\it spectral
decomposition\/} of a character. The measure $P$ will be called the
{\it spectral measure\/} of $\psi$. If $\psi$ is extreme then
its spectral measure reduces to the Dirac mass at $\psi$.

Let $(T,\xi)$ be a spherical representation of $(G,K)$, $\psi$ be the
corresponding character, and $P$ be its spectral measure. If $\xi$ is
replaced by another spherical vector in the same representation then
the character $\psi$ is changed, hence the measure $P$ is changed, too.
However, $P$ is transformed to an equivalent measure. Thus, the
equivalence class of $P$ is an invariant of $T$ as a unitary
representation.

The spectral decomposition of $\psi$ determines a decomposition of
the representation $T$ into a continual integral of irreducible
spherical representations.

\section{Thoma's theorem and spectral decomposition\\ of the
representations $T_z$ with $z\in\Z$}

\subsection{First example of extreme characters}

Let $\al=(\al_1\ge\dots\ge\al_p\ge0)$ and
$\be=(\be_1\ge\dots\ge\be_q\ge0)$ be two
collections of numbers such that
$$
\sum_{i=1}^p \al_i+\sum_{j=1}^q\be_j=1.
$$
Here one of the numbers $p,q$ may be zero (then the corresponding
collection $\al$ or $\be$ disappears). To these data we will assign an
extreme character $\chi^{(\al,\be)}$ of $S(\infty)$, as follows.

Let
$$
p_k(\al,\be)=\sum_{i=1}^p \al_i^k+(-1)^{k-1}\sum_{j=1}^q\be_j^k.
$$
Note that
$$
p_1(\al,\be)\equiv1.
$$
Given $\sigma\in S(\infty)$, we denote by $m_k(\sigma)$ the number of
$k$--cycles in $\sigma$. Since $\sigma$ is a finite permutation, we
have
$$
m_1(\sigma)=\infty, \qquad
\text{$m_k(\sigma)<\infty$ for $k\ge2$}, \qquad
\text{$m_k(\sigma)=0$ for $k$ large enough}.
$$
In this notation, we set
$$
\chi^{(\al,\be)}(\sigma)=\prod_{k=1}^\infty
(p_k(\al,\be))^{m_k(\sigma)}=\prod_{k=2}^\infty
(p_k(\al,\be))^{m_k(\sigma)}, \qquad \sigma\in S(\infty),
$$
where we agree that $1^\infty=1$ and $0^0=1$.

\begin{proposition*} Each function $\chi^{(\al,\be)}$ defined
by the above formula is an extreme character of $S(\infty)$.
\end{proposition*}

If $p=1$ and $q=0$ (i.e., $\al_1=1$ and all other parameters disappear) then we
get the trivial character, which equals 1 identically. If $p=0$ and $q=1$ then
we get the alternate character $\sgn(\sigma)=\pm1$, where the plus--minus sign
is chosen according to the parity of the permutation. More generally, we have
$$
\chi^{(\al,\be)}\cdot\sgn=\chi^{(\be,\al)}.
$$

\subsection{Thoma's set}

Let $\R^\infty$ denote the direct product of countably many copies of
$\R$. We equip $\R^\infty$ with the product topology. Let $\Om$ be
the subset of $\R^\infty\times\R^\infty$ formed by couples
$\al\in\R^\infty$, $\be\in\R^\infty$ such that
$$
\al=(\al_1\ge\al_2\ge\dots\ge0), \quad
\be=(\be_1\ge\be_2\ge\dots\ge0), \quad
\sum_{i=1}^\infty \al_i+\sum_{j=1}^\infty\be_j\le1.
$$

We call $\Om$ the {\it Thoma set.\/} We equip it with topology
induced from that of the space $\R^\infty\times\R^\infty$.
It is readily seen that $\Om$ is a compact space.

The couples $(\al,\be)$ that we dealt with in \S3.1 can be viewed as
elements of $\Om$.  The subset of such couples (with given $p,q$)
will be denoted by $\Om_{pq}$.

Note that each $\Om_{pq}$ is isomorphic to a
simplex of dimension $p+q-1$. As affine coordinates of the
simplex one can take the numbers
$$
\al_1-\al_2,\,\dots,\, \al_{p-1}-\al_p, \,\al_p, \,
\be_1-\be_2,\,\dots, \,\be_{q-1}-\be_q, \, \be_q
$$
but we will not use these coordinates.

\begin{proposition*} The union of the simplices $\Om_{pq}$ is
dense in $\Om$.
\end{proposition*}

For instance, the point
$(\underline0,\underline0)=(\al\equiv0,\be\equiv0)\in\Om$ can be
approximated by points of the simplices $\Om_{p0}$ as $p\to\infty$,
$$
(\underline0,\underline0)=\lim_{p\to\infty}
((\underbrace{1/p,\dots,1/p}_p), \underline0).
$$

\subsection{Description of extreme characters}

Now we extend by continuity the definition of \S3.1. For any
$k=2,3,\dots$ we define the function $p_k$ on $\Om$ as follows. If
$\om=(\al,\be)\in\Om$ then
$$
p_k(\om)=p_k(\al,\be)
=\sum_{i=1}^\infty \al_i^k+(-1)^{k-1}\sum_{j=1}^\infty\be_j^k.
$$
Note that $p_k$ is a continuous function on $\Om$. It should be
emphasized that the condition $k\ge2$ is necessary here: the similar
expression with $k=1$ (that is, the sum of all coordinates) {\it is
not\/} continuous.

Next, for any $\om=(\al,\be)\in\Om$ we set
$$
\chi^{(\om)}(\sigma)=\chi^{(\al,\be)}(\sigma)=\prod_{k=2}^\infty
(p_k(\al,\be))^{m_k(\sigma)}, \qquad \sigma\in S(\infty),
$$

\newtheorem*{thth}{Thoma's theorem}
\begin{thth}{\rm(i)} For any $\om\in\Om$ the function
$\chi^{(\om)}$ defined above is an extreme character of $S(\infty)$.

{\rm(ii)} Each extreme character is obtained in this way.

{\rm(iii)} Different points $\om\in\Om$ define different characters.
\end{thth}

In particular, the character $\chi^{(\underline0,\underline0)}$ is
the delta function at $e\in S(\infty)$. It corresponds to the
biregular representation defined in \S1.3.

Note that the topology of $\Om$ agrees with the topology of pointwise
convergence of characters on $S(\infty)$. This implies, in particular, that the
characters of \S3.1 are dense in the whole set of extreme characters with
respect to the topology of pointwise convergence.

\begin{corollary*} For any character $\psi$ of $S(\infty)$, its
spectral measure $P$ can be viewed as a probability measure on the
compact space $\Om$, and the integral representation of\/ {\rm \S2.5}
can be rewritten in the following form
$$
\psi(\sigma)=\int_{\Om}\chi^{(\om)}P(d\om), \qquad \sigma\in S(\infty).
$$
\end{corollary*}

\sbox{\olshsb}{$z$}
\subsection{Spectral decomposition for integral values of \usebox{\olshsb}}

Consider the generalized regular representations $T_z$ of the group
$G$ introduced in \S1.7.

\begin{thm*} Assume $z$ is an integer, $z=k\in\Z$.

{\rm(i)} The representation $T_k$ possesses $K$--invariant cyclic
vectors, i.e., it can be made a spherical representation.

{\rm(ii)} Let $\xi$ be any such vector, $\psi$ be the corresponding character,
and $P$ be its spectral measure on $\Om$. Then $P$ is supported by the subset
$$
\bigcup_{p,q\ge0, \, (p,q)\ne(0,0), \, p-q=k} \Om_{pq}
$$
and for any $\Om_{pq}$ entering this subset, the restriction of $P$
to $\Om_{pq}$ is equivalent to Lebesgue measure on the simplex
$\Om_{pq}$.
\end{thm*}

When $k\ne0$, the restriction $(p,q)\ne(0,0)$ is redundant because it
follows from the condition $p-q=k$.

The condition $p-q=k$ also implies
that the spectral measures corresponding to different integral values
of the parameter $z$ are mutually singular. This, in turn, implies
that the corresponding representations are {\it disjoint\/}, i.e.,
they do not have equivalent subrepresentations.

\section{The characters $\chi_z$}

\sbox{\olshsb}{$\chi_z$}
\subsection{Definition of \usebox{\olshsb} and its explicit expression}

Let $T_z$ be a generalized regular representation of $G$. Assume first
$z\ne0$. Recall that $T_z$ is realized in the Hilbert space
$L^2(\bS,\mu_t)$, where $t=|z|^2$. Let $\one$ denote the function on
$\bS$ identically equal to 1. It can be viewed as a vector of
$L^2(\bS,\mu_t)$. Since $\mu_t$ is $K$--invariant and the cocycle
$z^{c(x,g)}$ entering the construction of $T_z$ is trivial on $K$,
the vector $\one$ is a $K$--invariant vector. Consider the
corresponding matrix coefficient and pass to the corresponding
character (see \S2.4), which we denote by $\chi_z$. Thus,
$$
\chi_z(\sigma)=(T_z(\sigma,e)\one,\one), \qquad \sigma\in S(\infty).
$$

We aim to give a formula for $\chi_z$. To do this we will describe
the expansion of $\chi_z\mid_{S(n)}$ in irreducible characters of
$S(n)$ for any $n$. Recall that the irreducible representations of
$S(n)$ are parametrized by Young diagrams with $n$ boxes. Let $\Y_n$
be the set of these diagrams. For $\la\in\Y_n$ we denote by
$\chi^\la$ the corresponding irreducible character (the trace of the
irreducible representation of $S(n)$ indexed by $\la$). Let
$\dim\la=\chi^\la(e)$ be the dimension of this representation. In
combinatorial terms, $\dim\la$ is the number of standard Young
tableaux of shape $\la$. Note that for this number there exist closed
expressions. Below the notation $(i,j)\in\la$ means that the box on
the intersection of the $i$th row and the $j$th column belongs to $\la$.

\begin{thm*} For any $n=1,2,\dots$,
$$
\chi_z\mid_{S(n)}
=\sum_{\la\in\Y_n}\left(\frac{\prod\limits_{(i,j)\in\la}|z+j-i|^2}
{|z|^2\,(|z|^2+1)\dots(|z|^2+n-1)}\,\frac{\dim\la}{n!}\right) \chi^\la.
$$
\end{thm*}

Note that this formula also makes sense for $z=0$.

The next claim is a direct consequence of the formula.

\begin{proposition*} The function $\one$ is a cyclic vector for $T_z$
if and only if $z\notin\Z$.
\end{proposition*}

Thus, for nonintegral $z$, the couple $(T_z,\one)$ is a spherical
representation and the character $\chi_z$ entirely determines $T_z$.

Note that for $z=k\in\Z$, the cyclic span of $\one$ is a proper
subrepresentation that ``corresponds'' to a particular simplex
$\Om_{pq}$ (see \S3.4). Specifically,
$$
(p,q)=\begin{cases} (k,0), & \text{if $k>0$}\\ (0,|k|), & \text{if $k<0$}\\
(1,1), & \text{if $k=0$}. \end{cases}
$$
\sbox{\olshsb}{$z\leftrightarrow\bar z$}
\subsection{The symmetry \usebox{\olshsb}}

\begin{proposition*} For any $z$, the representations $T_z$ and
$T_{\bar z}$ are equivalent.
\end{proposition*}

\begin{proof} Indeed, if $z\in\R$ then there is nothing to prove. If
$z\notin\R$ then $\one$ is cyclic, so that the claim follows from the
fact that $\chi_z=\chi_{\bar z}$, which in turn is evident from
Theorem of $\S4.1$.
\end{proof}

Note that this is by no means evident from the construction of the
representations $T_z$.

\subsection{Disjointness}

Let $P_z$ be the spectral measure of the character $\chi_z$, see
\S3.3. When $z$ is integral, the measure $P_z$ lives on a simplex
$\Om_{pq}$, see \S4.1. Now we focus on the measures $P_z$ with
$z\notin\Z$.

\begin{thm*} {\rm(i)} Let $z\notin\Z$. Then all simplices
$\Om_{pq}$ are null sets with respect to the measure $P_z$.

{\rm(ii)} Let $z_1$ and $z_2$ be two complex number, both nonintegral,
$z_1\ne z_2$, and $z_1\ne\bar z_2$. Then the measures $P_{z_1}$ and
$P_{z_2}$ are mutually singular.
\end{thm*}

It follows that the generalized regular representations $T_z$ are
mutually disjoint, with the exception of the equivalence $T_z\sim
T_{\bar z}$.

\subsection{A nondegeneracy property}

\begin{proposition*} All measures $P_z$, $z\in\C$, are supported by
the subset
$$
\Om_0:=\left\{(\al,\be)\bigm| \sum\al_i+\sum\be_j=1\right\}.
$$
\end{proposition*}

On the contrary, the measure $P_\infty$ that corresponds to the
biregular representation $T_\infty$ is the Dirac measure at the point
$(\underline0,\underline0)$, which is outside $\Om_0$. This does not
contradicts the fact that the family $\{T_z\}$ is a deformation of
$T_\infty$, because $\Om_0$ is dense in $\Om$.

\section{Determinantal point processes}

\subsection{Point configurations}

Let $\mathfrak X$ be a locally compact separable topological space. By a {\it
point configuration\/} in $\mathfrak X$ we mean a locally finite collection $C$
of points of the space $\mathfrak X$. These points will also be called {\it
particles.\/} Here ``locally finite'' means that the intersection of $C$ with
any relatively compact subset is finite. Thus, $C$ is either finite or
countably infinite. Multiple particles in $C$ are, in principle, permitted but
all multiplicities must of course be finite. However, we will not really deal
with configurations containing multiple particles. Let us emphasize that the
particles in $C$ are unordered.

The set of all point configurations in $\mathfrak X$ will be denoted by
$\Conf(\mathfrak X)$.

\subsection{Definition of a point process}

A relatively compact Borel subset $A\subset\mathfrak X$ will be called a {\it
window.\/} Given a window $A$ and $C\in\Conf(\mathfrak X)$, let $\Cal N_A(C)$
be the cardinality of the intersection $A\cap C$ (with multiplicities counted).
Thus, $\Cal N_A$ is a function on $\Conf(\mathfrak X)$ taking values in $\Z_+$.
We equip $\Conf(\mathfrak X)$ with the Borel structure generated by the
functions of the form $\Cal N_A$.

By a measure on $\Conf(\mathfrak X)$ we will mean a Borel measure with respect
to this Borel structure.

By definition, a {\it point process\/} on $\mathfrak X$ is a probability
measure $\Cal P$ on the space $\Conf(\mathfrak X)$.

In practice, point processes often arise as follows. Assume we are given a
Borel space $Y$ and a map $\phi: Y\to\Conf(\mathfrak X)$. The map $\phi$ must
be a Borel map. i.e., for any window $A$, the superposition $\Cal N_A\circ\phi$
must be a Borel function on $Y$. Further, assume we are given a probability
Borel measure $P$ on $Y$. Then its pushforward $\Cal P$ under $\phi$ is well
defined and it is a point process.

Given a point process, we can speak about {\it random\/} point
configurations $C$. Any reasonable (that is, Borel) function of $C$
becomes a random variable. For instance, $\Cal N_A$ is a random
variable for any window $A$, and we may consider the probability that
$\Cal N_A$ takes any prescribed value.

\subsection{Example: Poisson process}

Let $\rho$ be a measure on $\mathfrak X$. It may be infinite but must take
finite values on any window. The {\it Poisson process with density\/} $\rho$ is
characterized by the following properties:

(i) For any window $A$, the random variable $\Cal N_A$ has the Poisson
distribution with parameter $\rho(A)$, i.e.,
$$
\Prob\{\Cal N_A=n\}=\frac{\rho(A)^n}{n!}e^{-\rho(A)}, \quad n\in\Z_+\,.
$$

(ii) For any pairwise disjoint windows $A_1,\dots,A_k$, the
corresponding random variables are independent.

In particular, if $\mathfrak X=\R$ and $\rho$ is the Lebesgue measure then this
is the classical Poisson process.

\subsection{Correlation measures and correlation functions}

Let $\Cal P$ be a point process on $\mathfrak X$. One can assign to $\Cal P$ a
sequence $\rho_1,\rho_2,\dots$ of measures, where, for any $n$, $\rho_n$ is a
symmetric measure on the $n$--fold product $\mathfrak X^n=\mathfrak
X\times\dots\times\mathfrak X$, called the {\it $n$--particle correlation
measure.\/} Under mild assumptions on $\Cal P$ the correlation measures exist
and determine $\Cal P$ uniquely. They are defined as follows.

Given $n$ and a compactly supported bounded Borel function $f$ on $\mathfrak
X^n$, let $\widetilde f$ be the function on $\Conf(\mathfrak X)$ defined by
$$
\widetilde f(C)=\sum_{i_1,\dots,i_n}f(x_{i_1},\dots,x_{i_n}), \qquad
C=\{x_1,x_2,\dots\}\in\Conf(\mathfrak X),
$$
summed over all $n$--tuples of pairwise distinct indices. Here we have
used an enumeration of the particles in $C$ but the result does not
depend on it.

Then the measure $\rho_n$ is characterized by the equality
$$
\int_{\mathfrak X^n}f\rho_n=\int_{C\in\Conf(\mathfrak X)}\widetilde f(C)\,\Cal
P(dC),
$$
where $f$ is an arbitrary compactly supported bounded Borel function on
$\mathfrak X^n$.

\newtheorem*{exams}{Examples}
\begin{exams} {\rm (i) If $\Cal P$ is a Poisson process then
$\rho_n=\rho^{\otimes n}$, where $\rho$ is the density of $\Cal P$.

(ii) Assume that $\mathfrak X$ is discrete and $\Cal P$ lives on multiplicity
free configurations. Then the correlation measures say what is the probability
that the random configuration contains an arbitrary given finite set of
points.}
\end{exams}

Often there is a natural measure $\nu$ on $\mathfrak X$ (a {\it reference
measure\/}) such that each $\rho_n$ has a density with respect to $\nu^{\otimes
n}$. This density is called the $n$th {\it correlation function.\/} For
instance, if $\mathfrak X$ is a domain of an Euclidean space and $\nu$ is the
Lebesgue measure then, informally, the $n$th correlation function equals the
density of the probability that the random configuration has particles in given
$n$ infinitisemal regions $dx_1$, \dots, $dx_n$.

\subsection{Determinantal point processes}

Let $\Cal P$ be a point process on $\mathfrak X$. Assume that $\mathfrak X$ is
equipped with a reference measure $\nu$ such that the correlation functions
(taken with respect to $\nu$) exist. Let us denote these functions by
$\rho_n(x_1,\dots,x_n)$. The process $\Cal P$ is said to be {\it
determinantal\/} if there exists a function $K(x,y)$ on $\mathfrak
X\times\mathfrak X$ such that
$$
\rho_n(x_1,\dots,x_n)=\det[K(x_i,x_j)]_{i,j=1}^n, \qquad
n=1,2,\dots\,.
$$
Then $K(x,y)$ is called a {\it correlation kernel\/} of $\Cal P$.

If $K(x,y)$ exists it is not unique since for any nonvanishing function
$\phi(x)$ on $\mathfrak X$, the kernel $\phi(x)K(x,y)(\phi(y))^{-1}$ leads to
the same result.

If we replace the reference measure by an equivalent one then we always can
appropriately change the kernel. Specifically, if $\nu$ is multiplied by a
positive function $f(x)$ then $K(x,y)$ can be replaced, say, by
$K(x,y)(f(x)f(y))^{-1/2}$.

\begin{exams} {\rm (i) Let $\mathfrak X=\R$, $\nu$ be the Lebesgue
measure, and $K(x,y)=\overline{K(y,x)}$ be the kernel of an Hermitian
integral operator $K$ in $L^2(\R)$. Then $K(x,y)$ is a correlation
kernel of a determinantal point process if and only if $0\le K\le1$ and
the restriction of the kernel to any bounded interval determines a
trace class operator.

(ii) The above conditions are satisfied by the {\it sine kernel}
$$
K(x,y)=\frac{\sin(\pi(x-y))}{\pi(x-y)}\,, \qquad x,y\in\R.
$$
The sine kernel arises in random matrix theory. It determines a translation
invariant point process on $\R$, which is a fundamental and probably the best
known example of a determinantal point process. }
\end{exams}

\section{The point processes $\Cal P_z$ and $\Pz$. The main result}
\nopagebreak

\subsection{From spectral measures to point processes}

Let $I=[-1,1]\subset\R$ and $I^*=[-1,1]\setminus\{0\}$. Let us take $I^*$ as
the space $\mathfrak X$. We define an embedding $\Om\to\Conf(I^*)$ as follows
$$
\om=(\al,\be)\mapsto C=\{\al_i\ne0\}\cup\{-\be_j\ne0\}.
$$
That is, we remove the possible zero coordinates, change the sign of the
$\beta$--coordinates, and forget the ordering.  In this way we convert $\om$ to
a point configuration $C$ in the punctured segment $I^*$. In particular, the
empty configuration $C=\varnothing$ corresponds to
$\om=(\underline0,\underline0)$.

Given a probability measure $P$ on $\Om$, its pushforward under this
embedding is a probability measure $\Cal P$ on $\Conf(I^*)$, i.e., a
point process on the space $I^*$, see \S5.2. Applying this procedure to
the spectral measures $P_z$ (\S4.3) we get point processes $\Cal
P_z$ on $I^*$.

\subsection{Lifting}

We aim to define a modification of the point processes $\Cal P_z$.
Fix $z\in\C\setminus\{0\}$ and set as usual $t=|z|^2$. Let $s>0$ be a
random variable whose distribution has the form
$$
\frac1{\Ga(t)}s^{t-1}e^{-s}ds
$$
(the gamma distribution on $\R_+$ with parameter $t$.) We assume that
$s$ is independent of $\Cal P_z$. Given the
random configuration $C$ of the process $\Cal P_z$, we multiply the
coordinates of all particles of $C$ by the random factor $s$. The
result is a random point configuration $\widetilde C$ on
$\R^*=\R\setminus\{0\}$.

We call this procedure the {\it lifting.\/} Under the lifting the
point process $\Cal P_z$ is transformed to a point process on $\R^*$
which we denote by $\Pz$.

The lifting is in principle reversible. Indeed, due to Proposition of
\S4.4, we can recover $C$ from $\widetilde C$ by dividing all the
coordinates in $\widetilde C$ by the sum of their absolute values.

It turns out that the lifting leads to a simplification of the
initial point process.

\subsection{Transformation of the correlation functions under the
lifting}

Fix the parameter $z$. Let $\rho_n(x_1,\dots,x_n)$ and
$\widetilde\rho_n(x_1,\dots,x_n)$ be the correlation
functions of the processes $\Cal P_z$ and $\Pz$, respectively (see
\S5.4).  Here we take the Lebesgue measure as the reference measure.

The definition of the lifting implies that
$$
\widetilde\rho_n(x_1,\dots,x_n)=
\int_0^\infty\frac{s^{t-1}e^{-s}}{\Ga(t)}
\,\rho_n\left(\frac{x_1}s,\dots,\frac{x_n}s\right)\,\frac{ds}{s^n}\,,
$$
where we agree that the function $\rho_n$ vanishes on
$(\R^*)^n\setminus(I^*)^n$. Thus, the action of the lifting on the
correlation functions is expressed by a ray integral transform.

This ray transform can be readily reduced to the Laplace transform.
It follows that it is injective, which agrees with the fact that
lifting is reversible.

\subsection{The main result}

To state the result we need some notation.

Let $W_{\kappa,\mu}(x)$ denote the {\it Whittaker function\/} with
parameters $\kappa,\mu\in\C$. It is a unique solution of the
differential equation
$$
W''-\left(\frac14-\frac{\kappa}{x}+ \frac{\mu^2-\frac14}{x^2}\right)W=0
$$
with the condition $W(x)\sim x^\kappa e^{-\frac{x}2}$ as $x\to+\infty$.

This function is initially defined for real positive
$x$ and then can be extended to a holomorphic function on
$\C\setminus(-\infty,0]$.

Next, we write $z=a+ib$ with real $a,b$ and set
$$
P_\pm(x)=\frac{t^{\frac12}}{|\Ga(1\pm z)|}\,
W_{\pm a+\frac12,ib}(x), \qquad
Q_\pm(x)=\frac{t^{\frac32}x^{-\frac12}}{|\Ga(1\pm z)|}\,
W_{\pm a-\frac12,ib}(x)\, .
$$

\newtheorem*{mnth}{Main Theorem}
\begin{mnth}For any $z\in\C\setminus\{0\}$, the point process
$\Pz$ is a determinantal process whose correlation kernel can be
written as
$$
K(x,y)=\begin{cases}\dfrac{P_+(x)Q_+(y)-Q_+(x)P_+(y)}{x-y}\,, &  x>0\,,y>0\\
\dfrac{P_+(x)P_-(-y)+Q_+(x)Q_-(-y)}{x-y}\,, & x>0\,,y<0\\
\dfrac{P_+(x)P_+(y)+Q_-(-x)Q_+(y)}{x-y}\,, & x<0\,,y>0\\
-\dfrac{P_-(-x)Q_-(-y)-Q_-(-x)P_-(-y)}{x-y}\,, & x<0\,,y<0\end{cases}
$$
where $x,y\in\R^*$ and the indeterminacy arising for $x=y$ is
resolved via the L'Hospital rule.
\end{mnth}

We call the kernel $K(x,y)$ the {\it Whittaker kernel.\/}

Note that $K(x,y)$ is real valued. It is not symmetric but satisfies
the symmetry property
$$
K(x,y)=\sgn(x)\sgn(y)K(y,x),
$$
where $\sgn(x)$ equals $\pm1$ according to the sign of $x$.
This property can be called {\it $J$--symmetry}, it means that
the kernel is symmetric with respect to an indefinite inner product.

\subsection{The $L$--operator}

Split the Hilbert space $L^2(\R^*)$ into the direct sum
$L^2(\R_+)\oplus L^2(\R_-)$, where all $L^2$ spaces are taken with
respect to the Lebesgue measure. According to this splitting we will
write operators in $L^2(\R^*)$ in block form, as $2\times2$ operator
matrices.

Let
$$
L=\bmatrix 0 & A\\ -A^t & 0 \endbmatrix\,,
$$
where $A$ is the integral operator with the kernel
$$
A(x,y)=\left|\frac{\sin(\pi z)}{\pi}\right|^2\cdot
\dfrac{\left(\frac{x}{|y|}\right)^{\operatorname{Re}z}
e^{-\frac{x-y}2}} {x-y}\,, \quad x>0, \quad y<0.
$$
By $A^t$ we denote the conjugate operator $L^2(\R_+)\to L^2(\R_-)$.

\begin{thm*} Assume that $-\frac12<\Re z<\frac12, \, z\ne0$.
Then $A$ is a bounded operator $L^2(\R_-)\to L^2(\R_+)$ and the
correlation kernel $K(x,y)$ is the kernel of the operator $L(1+L)^{-1}$.
\end{thm*}

Note that, in contrast to $K$, the kernel of $L$ does not involve
special functions.

\subsection{An application}

Fix $z\in\C\setminus\Z$ and consider the probability space
$(\Om,P_z)$. For any $k=1,2,\dots$ the coordinates $\al_k$
and $\be_k$ are functions in $\om\in\Om$, hence we may view them as
random variables. The next result provides an information about the
rate of their decay as $i, j\to\infty$.

\begin{thm*} With probability 1, there exist limits
$$
\lim_{k\to\infty}(\al_k)^{\frac 1k}
=\lim_{k\to\infty}(\be_k)^{\frac 1k}
=q(z)\in(0,1),
$$
where
$$
q(z)=\exp\left(\pi\,\frac{\ctg \pi z-\ctg\pi\bar z}{z-\bar z}\right)
=\exp\left(-\sum_{n\in\Z}\frac1{|z-n|^2}\right)
$$
\end{thm*}

\section{Scheme of the proof of the Main Theorem}

\sbox{\olshsb}{$z$}
\subsection{The z--measures}

Recall that by $\Y_n$ we denote the finite set of Young diagrams with
$n$ boxes. Set
$$
P_z^\n(\la)
=\frac{\prod\limits_{(i,j)\in\la}|z+j-i|^2}
{|z|^2\,(|z|^2+1)\dots(|z|^2+n-1)}\,\frac{(\dim\la)^2}{n!}\,,
\qquad \la\in\Y_n\,.
$$
Comparing this with the expression of $\chi_z\mid_{S(n)}$ (\S4.1) we see that
the quantities $P_z^\n(\la)$ are the coefficients in the expansion of $\chi_z$
in the {\it normalized\/} irreducible characters $\chi^\la/\dim\la$. It follows
that
$$
\sum_{\la\in\Y_n}P_z^\n(\la)=1.
$$
Thus, for any fixed $n=1,2,\dots$, the quantities $P_z^\n(\la)$
determine a probability measure on $\Y_n$. We will denote it by
$P_z^\n$ and call it the {\it z--measure\/} on $Y_n$.

\sbox{\olshsb}{$\Y_n\hookrightarrow\Om$}
\subsection{Frobenius coordinates and the embedding
\usebox{\olshsb}}

Given $\la\in\Y_n$, let $\la'$ be the transposed diagram and $d$ be
the number of diagonal boxes in $\la$.
We define the {\it modified Frobenius coordinates\/} of $\la$ as
$$
a_i=\la_i-i+\frac12,\qquad b_i=\la'_i-i+\frac12, \qquad i=1,\dots,d.
$$
Note that
$$
a_1>\dots>a_d>0, \qquad b_1>\dots>b_d>0, \qquad
\sum_{i=1}^d(a_i+b_i)=n.
$$

For any $n=1,2,\dots$ we embed $\Y_n$ into $\Om$ by making use of the
map
\begin{gather*}
\la\mapsto\om_\la=(\al,\be), \\
\al=\left(\frac{a_1}n,\,\dots,\,\frac{a_d}n,\,0,\,0,\,\dots\right), \qquad
\be=\left(\frac{b_1}n,\,\dots,\,\frac{b_d}n,\,0,\,0,\,\dots\right).
\end{gather*}

As $n\to\infty$, the points $\om_\la$ coming from the diagrams $\la\in\Y_n$
fill out the space $\Om$ more and more densely. Thus, for large $n$, the image
of $\Y_n$ in $\Om$ can be viewed as a discrete approximation of $\Om$.

\newsavebox{\olshsbii}\sbox{\olshsbii}{$z$}\sbox{\olshsb}{$P_z$}
\subsection{Approximation of \usebox{\olshsb} by z--measures}

Let $\underline P_z^\n$ be the pushforward of the measure $P_z^\n$
under the embedding $\Y_n\hookrightarrow\Om$. This is a probability
measure on $\Om$.

\newtheorem*{apprth}{Approximation Theorem}
\begin{apprth}As $n\to\infty$, the measures
$\underline P_z^\n$ weakly converge to the measure $P_z$.
\end{apprth}

This fact is the starting point for explicit computations related to
the measures $P_z$.

\subsection{The mixed z--measures}

Let $\Y=\Y_0\cup\Y_1\cup\Y_2\cup\dots$ be the set of all Young diagrams.
We agree that $\Y_0$ consists of a single element -- the empty
diagram $\varnothing$. Fix $z\in\C\setminus\{0\}$ and $\xi\in(0,1)$.
We define a measure $\widetilde P_{z,\xi}$ on $\Y$ as
follows:
$$
\widetilde P_{z,\xi}(\la)=P_z^\n(\la)\cdot
(1-\xi)^{|z|^2}\, \frac{|z|^2\,(|z|^2+1)\dots(|z|^2+n-1)}{n!}\,\xi^n,
\qquad \la\in\Y,
$$
where $n$ is the number of boxes in $\la$ and
$P_z^{(0)}(\varnothing):=1$.

In other words, $\widetilde P_{z,\xi}$ is obtained by mixing together all the
z--measures $P^{(0)}_z,P^{(1)}_z,\dots$, where the weight of the $n$th
component is equal to
$$
\pi_{t,\xi}(n)=(1-\xi)^t\, \frac{t(t+1)\dots(t+n-1)}{n!}\,\xi^n,
\qquad t=|z|^2.
$$
Note that
$$
\sum_{n=0}^\infty \pi_{t,\xi}(n)=1.
$$
It follows that $\widetilde P_{z,\xi}$ is a probability measure. Let
us call it the {\it mixed z--measure.\/}

Note that, as $z\to0$, the measure $\widetilde P_{z,\xi}$ tends to
the Dirac mass at $\{\varnothing\}$ for any fixed $\xi$.
\sbox{\olshsb}{$\widetilde{\Cal P}_{z,\xi}$}
\subsection{The lattice process \usebox{\olshsb}}

Set
$$
\Z'=\Z+\frac12
=\left\{\dots,\,-\frac32,\,-\frac12,\,\frac12,\,\frac32,\,\dots\right\}.
$$
Using the notation of \S7.2 we assign to an arbitrary Young diagram a point
configuration $C\in\Conf(\Z')$, as follows
$$
\la\mapsto C=\{-b_1,\dots,-b_d,a_d,\dots,a_1\}.
$$
The correspondence $\la\mapsto C$ defines an embedding
$\Y\hookrightarrow\Conf(\Z')$. Take the pushforward of the measure
$P_{z,\xi}$ under this embedding. It is a probability measure on
$\Conf(\Z')$, hence a point process on $\Z'$. Let us denote it by
$\widetilde{\Cal P}_{z,\xi}$.

\begin{thm*} The process $\widetilde{\Cal P}_{z,\xi}$ on the
lattice $\Z'$ is determinantal. Its correlation kernel can be
explicitly computed: it has the form quite similar to that of the kernel
$K(x,y)$ from \S6.4, where the corresponding functions $P_\pm$ and
$Q_\pm$ are now expressed through the Gauss hypergeometric function.
\end{thm*}

\subsection{Idea of proof of the Main Theorem}

Given $\xi\in(0,1)$, we embed the lattice $\Z'$ into $\R^*$ as follows
$$
\Z'\ni x\mapsto(1-\xi)x\in\R^*.
$$
Let $\widetilde{\underline{\Cal P}}_{z,\xi}$ be the pushforward of
$\widetilde{\Cal P}_{z,\xi}$ under this embedding. We can view
$\widetilde{\underline{\Cal P}}_{z,\xi}$ as a point process on
$\R^*$.

Remark that the probability distribution $\pi_{t,\xi}$ on $\Z_+$
introduced in \S7.4 approximates in an appropriate scaling limit as
$\xi\nearrow1$ the gamma distribution on $\R^*$ with parameter $t$.
Specifically, the scaling has the form $n\mapsto (1-\xi)n$.
Recall that we have used the gamma distribution in the definition of
the lifting, see \S6.2.

Combining this fact with the Approximation Theorem of \S7.3 we
conclude that the process $\widetilde{\underline{\Cal P}}_{z,\xi}$
must converge to the process $\Pz$ as $\xi\nearrow1$ in a certain
sense.  More precisely, we prove that the correlation measures of the
former process converge to the respective correlation measures of the
latter process.

On the other hand, we can explicitly compute the scaled limit of the
lattice correlation measures using the explicit expression of the
lattice kernel from \S7.5. It turns out that then the Gauss
hypergeometric function degenerates to the Whittaker function and we
get the formulas of \S6.4.

\section{Notes and references}

\subsection{Section 1}

The main reference to this section is the paper
Kerov--Olshanski--Vershik \cite{ol-KOV}.

\S1.1. Peter--Weyl's theorem is included in many textbooks on
representation theory. See, e.g., Naimark \cite{ol-Na},~\S32.

\S1.2. From the purely algebraic point of view, there is no single infinite
analog of the permutation groups $S(n)$ but a number of different versions. The
group $S(\infty)=\varinjlim S(n)$ formed by {\it finite\/} permutations of the
set $\{1,2,\dots\}$ and the group of {\it all\/} permutations of this set may
be viewed as the minimal and the maximal versions. There is also a huge family
of intermediate groups. The choice of an appropriate version may vary depending
on the applications we have in mind. Certain {\it topological\/} groups
connected with $S(\infty)$ are discussed in Olshanski \cite{ol-Ol1}, Okounkov
\cite{ol-Ok1}, \cite{ol-Ok1a}.

\S1.3. The result of the Proposition is closely related to von Neumann's
classical construction of II$_1$ factors. See Murray--von Neumann
\cite{ol-MvN}, ch. 5, and  Naimark \cite{ol-Na}, ch.~VII, \S38.5.

\S1.4. The $G$--space $\mathfrak S$ of virtual permutations was introduced in
Kerov--Olshanski--Vershik \cite{ol-KOV}. Notice that the canonical projection
$p_n$ emerged earlier, see Aldous \cite{ol-Ald}, p. 92. A closely related
construction, which also appeared earlier, is the so--called {\it Chinese
restaurant process,\/} see, e.g., Arratia--Barbour--Tavar\'e \cite{ol-ABD1},
\S2 and references therein. Projective limit constructions for classical groups
and symmetric spaces are considered in Pickrell \cite{ol-Pi}, Neretin
\cite{ol-Ne}, Olshanski \cite{ol-Ol4}. Earlier papers: Hida--Nomoto
\cite{ol-HN}, Yamasaki \cite{ol-Ya1}, \cite{ol-Ya2}, Shimomura \cite{ol-Shi}.

\S1.5. The definition of the Ewens measures $\mu_t$ on the space $\bS$ was
proposed in \cite{ol-KOV}, see also Kerov--Tsilevich \cite{ol-KT}, Kerov
\cite{ol-Ke4} (the latter paper deals with a generalization of these measures).
The definition of \cite{ol-KOV} was inspired by the fundamental concept of the
{\it Ewens sampling formula\/}, which was derived in 1972 by Ewens
\cite{ol-Ew1} in the context of population genetics. There is a large
literature concerning Ewens' sampling formula (or Ewens' partition structure).
See, e.g., the papers Watterson \cite{ol-Wat}, Kingman \cite{ol-Ki2},
\cite{ol-Ki3}, \cite{ol-Ki5}, Arratia--Barbour--Tavar\'e
\cite{ol-ABD1},\cite{ol-ABD2}, Ewens \cite{ol-Ew2}, which contain many other
references.

\S1.6. The results were established in \cite{ol-KOV}. For projective
limits of classical groups and symmetric spaces, there also
exist distinguished families of measures with good transformation
properties, see Pickrell \cite{ol-Pi}, Neretin \cite{ol-Ne},
Olshanski \cite{ol-Ol4}.

\S1.7. The representations $T_z$ were introduced in \cite{ol-KOV}. A parallel
construction exists for infinite--dimensional classical groups and symmetric
spaces, see the pioneer paper Pickrell \cite{ol-Pi} and also Neretin
\cite{ol-Ne}, Olshanski \cite{ol-Ol4}.

\subsection{Section 2}

\S2.1. The concept of spherical representations is usually employed for {\it
Gel\-fand pairs\/} $(\Cal G, \Cal K)$. According to the conventional
definition, $(\Cal G,\Cal K)$ is said to be a Gelfand pair if the subalgebra of
$\Cal K$--biinvariant functions in the group algebra $L^1(\Cal G)$ is
commutative. This works for locally compact $\Cal G$ and compact $\Cal K$.
There exists, however, a reformulation which makes sense for arbitrary groups,
see Olshanski \cite{ol-Ol2}. Our pair $(G,K)$ is a Gelfand pair, see Olshanski
\cite{ol-Ol3}.

\S2.2. For general facts concerning positive definite functions on
groups, see, e.g., Naimark \cite{ol-Na}.

\S2.3. There exist at least two different ways to define characters for
infinite--dimensional representations. The most known recipe (Gelfand,
Harish--Chandra) is to view characters not as ordinary functions but as
distributions on the group. This idea works perfectly for a large class of Lie
groups and $p$--adic groups but not for groups like $S(\infty)$. The definition
employed here follows another approach, which goes back to von Neumann. Extreme
characters of a group $\Cal K$ are related to finite factor representations of
$\Cal K$ in the sense of von Neumann. See Thoma \cite{ol-Th1}, \cite{ol-Th2},
Stratila--Voiculescu \cite{ol-SV}, Voiculescu \cite{ol-Vo}.

\S2.4. The correspondence between extreme characters and irreducible
spherical representations was pointed out in Olshanski \cite{ol-Ol1},
\cite{ol-Ol2}. The Proposition follows from the fact that our pair $(G,K)$
is a Gelfand pair, see Olshanski \cite{ol-Ol3}.

The irreducible spherical representations of $(G,K)$ form a subfamily of a
larger family of representations called {\it admissible representations,\/} see
Olshanski \cite{ol-Ol1}, \cite{ol-Ol2}, \cite{ol-Ol3}. On the other hand, aside
from {\it finite\/} factor representations of $S(\infty)$ that correspond to
extreme characters, there exist interesting examples of factor representations
of quite different nature, see Stratila--Voiculescu \cite{ol-SV}.

Explicit realizations of finite factor representations of $S(\infty)$
and irreducible spherical representations of $(G,K)$ are given in
Vershik--Kerov \cite{ol-VK1}, Wassermann \cite{ol-Was}, Olshanski \cite{ol-Ol3}.

\S2.5. There are various methods to establish the existence and uniqueness of
the spectral decomposition. See, e.g., Diaconis--Freedman \cite{ol-DF},
Voiculescu \cite{ol-Vo}, Olshanski \cite{ol-Ol4}. One more approach, which is
specially adapted to the group $S(\infty)$ and provides an explicit description
of $\Ex(\Cal X(S(\infty)))$, is proposed in Kerov--Okounkov--Olshanski
\cite{ol-KOO}.

\subsection{Section 3}

\S3.1. The expressions $p_k(\al,\be)$ are {\it supersymmetric\/}
analogs of power sums. About the role of supersymmetric functions in
the theory of characters of $S(\infty)$ see Vershik--Kerov \cite{ol-VK2},
Olshanski--Regev--Vershik \cite{ol-ORV}.

\S3.2. The Thoma set $\Om$ can be viewed as an infinite--dimensional
simplex. The subsets $\Om_{pq}$ are exactly its finite--dimensional
faces.

\S3.3. Thoma's paper \cite{ol-Th1} was the first work about characters of
$S(\infty)$. It contains the classification of extreme characters
(Thoma's theorem), which was obtained using complex--analytic tools.
Thoma's theorem is equivalent to another classification problem --- that
of {\it one--sided totally positive sequences.\/} Much earlier, that
problem was raised by Schoenberg and solved by Edrei \cite{ol-Ed}. The
equivalence of both problems was implicit in Thoma's paper \cite{ol-Th1}
but Thoma apparently was not aware of the works on total positivity.

The next step was made by Vershik and Kerov \cite{ol-VK2}. Following a general
principle earlier suggested in Vershik \cite{ol-Ve1}, Vershik and Kerov found a
new proof of Thoma's theorem. Their approach is based on studying the limit
transition from characters of $S(n)$ to characters of $S(\infty)$. This
provides a very natural interpretation of Thoma's parameters $\al_i$, $\be_j$.

Developing further the asymptotic approach of \cite{ol-VK2},
Kerov--Okounkov--Olshanski \cite{ol-KOO} obtained a generalization of Thoma's
theorem. An even more general claim was conjectured by Kerov in \cite{ol-Ke2}.

One of the fruitful ideas contained in Vershik--Kerov's paper \cite{ol-VK2}
concerns the combinatorics of irreducible characters $\;\chi^\la\;$ of the
finite symmetric groups. Assume that $\la\in\Y_n$ and $\rho$ is a partition of
$m$, where $m\le n$. Let $\chi^\la_\rho$ denote the value of $\chi^\la$ at the
conjugacy class in $S(n)$ indexed by the partition $\rho\cup1^{n-m}$ of $n$.
The idea was to consider $\chi^\la_\rho$ as a function in $\la$ with $\rho$
viewed as a parameter. Vershik and Kerov discovered that the function
$\la\mapsto\chi^\la_\rho$, after a simple normalization, becomes a
supersymmetric function in the modified Frobenius coordinates of $\la$. This
function is inhomogeneous and its top degree homogeneous term is the
supersymmetric (product) power sum indexed by $\rho$. Further results in this
directions: Kerov--Olshanski \cite{ol-KO}, Okounkov--Olshanski \cite{ol-OO},
Olshanski--Regev--Vershik \cite{ol-ORV}. Even in the simplest case when $\rho$
consists of a single part ($\rho=(m)$) the function
$\la\mapsto\chi^\la_\rho=\chi^\la_{(m)}$ is rather nontrivial. See Wassermann
\cite{ol-Was}, Kerov \cite{ol-Ke3}, Biane \cite{ol-Bi}, Ivanov--Olshanski
\cite{ol-IO}.

\S3.4. The spectral decomposition of $T_z$'s for integral values of $z$ was
obtained in Kerov--Olshanski--Vershik \cite{ol-KOV}.

\subsection{Section 4}

\S4.1. The results were obtained in Kerov--Olshanski--Vershik
\cite{ol-KOV}. Similar results for other groups: Pickrell \cite{ol-Pi},
Olshanski \cite{ol-Ol4}.

\S4.2. One can define intertwining operators for the representations
$T_z$ and $T_{\bar z}$. These operators have interesting properties.
See Kerov--Olshanski--Vershik \cite{ol-KOV}.

\S4.3. The result was obtained in Kerov--Olshanski--Vershik \cite{ol-KOV}.
Note that the Theorem of \S6.6 implies a weaker result: the spectral
measures $P_{z_1}$ and $P_{z_2}$ are mutually singular for any
$z_1,z_2\in\C\setminus\Z$ such that $q(z_1)\ne q(z_2)$.

\S4.4. The result was announced in Kerov--Olshanski--Vershik
\cite{ol-KOV}. It can be proved in different ways, see Olshanski
\cite{ol-P.I}, Borodin \cite{ol-P.II}.

\subsection{Section 5}

\S\S5.1 -- 5.4. The material is standard. See Daley and Vere-Jones
\cite{ol-DVJ}, Lenard \cite{ol-Len}, Kingman \cite{ol-Ki5}. Point processes are
also called {\it random point fields.}

\S5.5. The class of determinantal point process was first singled out by Macchi
\cite{ol-Ma1}, \cite{ol-Ma2} under the name of {\it fermion processes.\/} The
motivation comes from a connection with the fermionic Fock space. The term
``determinantal'' was suggested in Borodin--Olshanski \cite{ol-BO2}.  We found
it more appropriate, because in our concrete situation, point configurations
may be viewed as consisting of particles of {\it two opposite charges.\/} A
number of important examples of determinantal point processes emerged in random
matrix theory, see, e.g., Dyson \cite{ol-Dy}, Mehta \cite{ol-Me}, Nagao--Wadati
\cite{ol-NW}, Tracy--Widom \cite{ol-TW}, and the references therein.  However,
to our knowledge, up to the recent survey paper by Soshnikov \cite{ol-So}, the
experts in this field did not pay attention to general properties of
determinantal processes and did not introduce any general name for them.

The result stated in Example (i) is due to Soshnikov \cite{ol-So}.

\subsection{Section 6}

\S6.1. The spectral measures $P_z$ with nonintegral parameter $z$
originally looked mysterious: it was unclear how to handle them.

The idea of converting the measures $P_z$ into point
processes $\Cal P_z$ and computing the correlation functions was
motivated by the following observation. It turns out that the
coefficients of the expansion of \S4.1 can be interpreted as moments
of certain auxiliary measures (we called them the {\it controlling
measures\/}). The controlling measures are determined by these
moments uniquely. On the other hand, the correlation functions can be
expressed through the controlling measures. It follows that
evaluating the correlation functions can be reduced to solving
certain (rather complicated) multidimensional moment problems.

We followed first this way (see the preprints \cite{ol-P.I}--\cite{ol-P.V}; part of
results was published in Borodin \cite{ol-Bor2}, \cite{ol-Bor3}; a summary
is given in Borodin--Olshanski \cite{ol-BO1}). A general description of the
method and the evaluation of the first correlation function are given
in Olshanski \cite{ol-P.I}. In Borodin \cite{ol-Bor2} the moment problem in
question is studied in detail. This leads (Borodin \cite{ol-P.II}) to
some formulas for the higher correlation functions: a
multidimensional integral representation and an explicit expression
through a multivariate Laurichella hypergeometric series of type B.
Both are rather involved.

\S\S6.2--6.3. The idea of lifting (Borodin \cite{ol-P.IV}) turned out to be
extremely successful, because it leads to a drastic simplification of the
correlation functions. What is even more important is that due to this
procedure we finally hit a nice class of point processes, the determinantal
ones.

\S6.4. The derivation of the Whittaker kernel by the first method is given in
Borodin \cite{ol-P.IV}, \cite{ol-Bor3}. It should be noted that the Whittaker
kernel belongs to the class of {\it integrable kernels.\/} This class was
singled out by Its--Izergin--Korepin--Slavnov \cite{ol-IIKS}, see also Deift
\cite{ol-De}, Borodin \cite{ol-Bor6}.

\S6.5. The claim concerning the L--operator and some related facts
are contained in Olshanski \cite{ol-P.V}. A conclusion is that (at least
when $|\Re z|<1/2$) the whole information about the spectral measure
$P_z$ is encoded in a very simple kernel $L(x,y)$.

\S6.6. The result is obtained in Borodin--Olshanski \cite{ol-P.III}. It
can be viewed as a strong law of large numbers. Roughly speaking, the
coordinates $\al_k$, $\be_k$ decay like the terms of the geometric
progression $\{q(z)^k\}$. A similar result holds for point processes
of quite different type (Poisson--Dirichlet distributions), see
Vershik--Shmidt \cite{ol-VS}.

Notice that the preprints \cite{ol-P.I}--\cite{ol-P.V} contain a number of other
results, some of them remain still unpublished.

\subsection{Section 7}

The main reference for this section is the paper Borodin--Olshanski
\cite{ol-BO2}, which gives an alternate way of proving the Main Theorem.
The method of \cite{ol-BO2} is simpler than the previous approach based on
a moment problem. Furthermore, our second approach explains the
origin of the lifting. However, the correlation functions for the initial
process $\Cal P_z$ are not directly obtained in this way.

\S7.1. The z--measures $P_z^\n$ with fixed parameter $z$ and varying
index $n$ satisfy the {\it coherency relation\/}
$$
P_z^\n(\mu)=\sum_{\la\in\Y_{n+1}:\, \la\supset\mu}
\frac{\dim\mu}{\dim\la}\, P_z^{(n+1)}(\la), \qquad
n=1,2,\dots, \quad \mu\in\Y_n\,.
$$
It expresses the fact that the function $\chi_z\mid_{S(n+1)}$ is an extension
of the function $\chi_z\mid_{S(n)}$. The coherency relation is not evident from
the explicit expression for the z--measures.

As $|z|\to\infty$, the measures $P_z^\n$ converge to the {\it Plancherel
measure\/} on $\Y_n$,
$$
P_\infty^\n(\la)=\frac{(\dim\la)^2}{n!}\,.
$$
Note that the expression for $P_z^\n(\la)$ looks as  a
product over the boxes of $\la$ times $P_\infty^\n(\la)$. This
property together with the coherency relation can be used for a
combinatorial characterization of the z--measures, see Rozhkovskaya
\cite{ol-R}.

Actually, the term ``z--measures'' has a somewhat wider meaning: the
family $\{P_z^\n\}$ forms the ``principal series'' while the whole
family of the z--measures also includes a ``complementary series'' and a
``degenerate series'' of measures which are given by similar
expressions.

A much larger family of {\it Schur measures\/} was introduced by Okounkov
\cite{ol-Ok3}.  In general, the Schur measures do not obey the coherency
relation and hence do not correspond to characters of $S(\infty)$. However,
they also give rise to determinantal point processes. It would be interesting
to know whether the z--measures exhaust all Schur measures satisfying the
coherency relation.

Kerov \cite{ol-Ke5} introduced analogs of z--measures satisfying a certain
one--parameter {\it deformation\/} of the coherency relation (the coherency
relation written above is closely related to the Schur functions, while Kerov's
more general form of the coherency relation is related to the Jack symmetric
functions, see also Kerov--Okounkov--Olshanski \cite{ol-KOO}). For another
approach, see Borodin--Olshanski \cite{ol-BO3}. Study of the point processes
corresponding to these more general z--measures was started in
Borodin--Olshanski \cite{ol-BO7}.

An analog of z--measures corresponding to {\it projective\/} characters of
$S(\infty)$ was found in Borodin \cite{ol-Bor1}. See also Borodin--Olshanski
\cite{ol-BO3}.

The paper Borodin--Olshanski \cite{ol-BO4} presents a survey of connections
between z--measures and a number of models arising in combinatorics, tiling,
directed percolation and random matrix theory.

\S7.2. The idea of embedding $\Y$ into $\Om$ is due to Vershik and Kerov
\cite{ol-VK2}. In a more general context it is used in
Kerov--Okounkov--Olshanski \cite{ol-KOO}.

\S7.3. The Approximation Theorem actually holds for spectral measures
corresponding to arbitrary characters of $S(\infty)$. See
Kerov--Okounkov--Olshanski \cite{ol-KOO}.

\S7.4. What we called ``mixing'' is a well--known trick. Under different names
it is used in various asymptotic problems of combinatorics and statistical
physics. See, e.g., Vershik \cite{ol-Ve3}. The general idea is to replace a
large $n$ limit, where the index $n$ enumerates different probabilistic
ensembles, by a limit transition of another kind (we are dealing with a
unifying ensemble depending on a parameter and let the parameter tend to a
limit). In many situations the two limit transitions lead to the same result.
For instance, this usually happens for the {\it poissonization procedure,\/}
when the mixing distribution on $\Z_+$ is a Poisson distribution. (About the
poissonized Plancherel measure, see Baik--Deift--Johansson \cite{ol-BDJ},
Borodin--Okounkov--Olshanski \cite{ol-BOO}, Johansson \cite{ol-Jo}.) A key
property of the Poisson distribution is that as its parameter goes to infinity,
the standard deviation grows more slowly than the mean. In our situation,
instead of Poisson we have to deal with the distribution $\pi_{t,\xi}$, a
particular case of the {\it negative binomial distribution.\/} As
$\xi\nearrow1$, the standard deviation and the mean of $\pi_{t,\xi}$ have the
same order of growth, which results in a nontrivial transformation of the large
$n$ limit (the lifting).

\S7.5. The fact that the lattice process $\widetilde{\Cal P}_z$ is
determinantal is checked rather easily. The difficult part
of the Theorem is the calculation of the correlation kernel. This can
be done in different ways, see Borodin--Olshanski \cite{ol-BO2}, Okounkov
\cite{ol-Ok2}, \cite{ol-Ok3}. Borodin \cite{ol-Bor4}, \cite{ol-Bor6} describes a
rather general procedure of computing correlation kernels via a
Riemann--Hilbert problem.

\S7.6. For more details see Borodin--Olshanski \cite{ol-BO2}.

\subsection{Other problems of harmonic analysis leading to point
processes}

A parallel but more complicated theory holds for the infinite--dimensional
unitary group $U(\infty)=\varinjlim U(N)$. For this group, there exists a
completion $\mathfrak U$ of the group space $U(\infty)$, which plays the role
of the space $\bS$ of virtual permutations.  On $\mathfrak U$, there exists a
family of measures with good transformation properties which give rise to
certain unitary representations of $U(\infty)\times U(\infty)$ --- analogs of
the representations $T_z$. See Neretin \cite{ol-Ne}, Olshanski \cite{ol-Ol4}.
The problem of harmonic analysis for these representations is studied in
Borodin--Olshanski \cite{ol-BO6}. It leads to determinantal point processes on
the space $\R\setminus\{\pm\frac12\}$. Their correlation kernels were found in
\cite{ol-BO6}: these are integrable kernels expressed through the Gauss
hypergeometric function.

There exists a similarity between decomposition of unitary
representations into irreducible ones and decomposition of invariant
measures on ergodic components. Both problems often can be
interpreted in terms of barycentric decomposition on extreme points
in a convex set. Below we briefly discuss two problems of ``harmonic
analysis for invariant measures'' that lead to point processes.

The first problem concerns invariant probability measures for the action of the
diagonal group $K\subset G$ on the space $\bS$. Recall that $K$ is isomorphic
to $S(\infty)$. Such measures are in 1--1 correspondence with {\it partition
structures\/} in the sense of Kingman \cite{ol-Ki1}. The set of all
$K$--invariant probability measures on $\bS$ (or partition structures) is a
convex set. Its extreme points correspond to {\it ergodic invariant measures\/}
whose complete classification is due to Kingman \cite{ol-Ki3}, \cite{ol-Ki4},
see also Kerov \cite{ol-Ke1}. Kingman's result is similar to Thoma's theorem.
The decomposition of Ewens' measures $\mu_t$ on ergodic components leads to a
remarkable one--parameter family of point processes on $(0,1]$ known as {\it
Poisson--Dirichlet distributions.\/} There is a large literature on
Poisson--Dirichlet distributions, we cite only a few works: Watterson
\cite{ol-Wat}, Griffiths \cite{ol-Gri}, Vershik--Shmidt \cite{ol-VS}, Ignatov
\cite{ol-Ig}, Kingman \cite{ol-Ki1}, \cite{ol-Ki2}, \cite{ol-Ki5}, Vershik
\cite{ol-Ve2}, Arratia--Barbour--Tavar\'e \cite{ol-ABD2}. One can show that the
lifting of the Poisson--Dirichlet distribution with parameter $t>0$ is the
Poisson process on $(0,+\infty)$ with density $\frac txe^{-x}dx$.

In the second problem, one deals with $(U(\infty),\mathfrak U)$ instead of
$(S(\infty),\bS)$. Here we again have a distinguished family of invariant
measures, see Borodin--Olshanski \cite{ol-BO5}, Olshanski \cite{ol-Ol4}. Their
decomposition on ergodic components is described in terms of certain
determinantal point processes on $\R^*$. The corresponding correlation kernels
are integrable and are expressed through another solution of Whittaker's
differential equation (\S6.4), see \cite{ol-BO5}. This subject is closely
connected with Dyson's {\it unitary circular ensemble\/}, see \cite{ol-BO5},
\cite{ol-Ol4}.

For the point processes mentioned above, a very interesting quantity is the
position of the rightmost particle in the random point configuration. In the
Poisson--Dirichlet case, the distribution of this random variable is given by a
curious piece--wise analytic function satisfying a linear
difference--differential equation: see Vershik--Shmidt \cite{ol-VS}, Watterson
\cite{ol-Wat}. For the (discrete and continuous) determinantal point processes
arising in harmonic analysis, the distribution of the rightmost particle can be
expressed through solutions of certain nonlinear (difference or differential)
Painlev\'e equations: see Borodin \cite{ol-Bor5}, Borodin--Deift \cite{ol-BD}.


\begin{thebibliography}{99}

\bibitem{ol-Ald}
D.~J.~Aldous, \emph{Exchangeability and related topics,} in: Springer Lecture
Notes in Math. {\bf 1117} (1985), pp. 2--199.


\bibitem{ol-ABD1}
 R.~Arratia, A.~D.~Barbour, and S.~Tavar\'e,
\emph{Poisson processes approximations for the Ewens sampling formula,} Ann.
Appl. Probab., \textbf{2} (1992), 519--535.



\bibitem{ol-ABD2}
 R.~Arratia, A.~D.~Barbour, and S.~Tavar\'e,
\emph{Random combinatorial structures and prime factorizations,} Notices Amer.
Math. Soc., \textbf{44} (1997), no.~8, 903--910.


\bibitem{ol-BDJ}
J.~Baik, P.~Deift and K.~Johansson, \emph{On the distribution of the length of
the longest increasing subsequence of random permutations,} J. Amer. Math.
Soc., \textbf{12} (1999), no.~4, 1119--1178.



\bibitem{ol-Bi}
Ph. Biane, \emph{Representations of symmetric groups and free probability,}
Advances in Math., \textbf{138} (1998), 126--181.



\bibitem{ol-Bor1}
A.~Borodin, \emph{Multiplicative central measures on the Schur graph,} in:
Representation theory, dynamical systems, combinatorial and algorithmic methods
II (A.~M.~Vershik, ed.), Zapiski Nauchnykh Seminarov POMI \textbf{240}, Nauka,
St.~Petersburg, 1997, 44--52 (Russian); English translation: J. Math. Sci. (New
York), \textbf{96} (1999), no. 5, 3472--3477.



\bibitem{ol-Bor2}
A.~Borodin, \emph{Characters of symmetric groups and correlation functions of
point processes,} Funktsional. Anal. Prilozhen., \textbf{34} (2000), no.~1,
12--28 (Russian); English translation: Funct. Anal. Appl. \textbf{34} (2000),
no.~1, 10--23.



\bibitem{ol-Bor3}
A.~Borodin, \emph{Harmonic analysis on the infinite symmetric group and the
Whittaker kernel,} Algebra Anal. \textbf{12} (2001), no.~5, 28--63 (Russian);
English translation: St.~Petersburg Math. J., \textbf{12} (2001), no.~5,
733--759.



\bibitem{ol-Bor4}
A.~Borodin, \emph{Riemann--Hilbert problem and the discrete Bessel kernel,}
Intern. Math. Research Notices, \textbf{2000:9} (2000), 467--494; {\tt
math/9912093}.



\bibitem{ol-Bor5}
A.~Borodin, \emph{Discrete gap probabilities and discrete Painlev\'e
equations,} Duke Math. J, to appear;  {\tt math-ph/0111008}.



\bibitem{ol-Bor6}
A.~Borodin, \emph{Asymptotic representation theory and Riemann--Hilbert
problem,}  in this volume; {\tt math/0110318}.


\bibitem{ol-BD}
A.~Borodin and P.~Deift, \emph{Fredholm determinants, Jimbo--Miwa--Ueno
tau--functions, and representation theory,} Comm. Pure Appl. Math. \textbf{55}
(2002), no.~9, 1160-1230; {\tt math-ph/0111007}.



\bibitem{ol-BOO}
A.~Borodin, A.~Okounkov and G.~Olshanski, \emph{Asymptotics of Plancherel
measures for symmetric groups,} J. Amer. Math. Soc., \textbf{13} (2000), no.~3,
481--515; {\tt math/9905032}.



\bibitem{ol-BO1}
A.~Borodin and G.~Olshanski, \emph{Point processes and the infinite symmetric
group,} Math. Research Lett., \textbf{5} (1998), 799--816; {\tt math/9810015}.



\bibitem{ol-BO2}
A.~Borodin and G.~Olshanski, \emph{Distributions on partitions, point
processes, and the hypergeometric kernel,} Commun. Math. Phys., \textbf{211}
(2000), 335--358; {\tt math/9904010}.



\bibitem{ol-BO3}
A.~Borodin and G.~Olshanski, \emph{Harmonic functions on multiplicative graphs
and interpolation polynomials,} Electronic J. Comb., \textbf{7} (2000), paper
\#R28; {\tt math/9912124}.



\bibitem{ol-BO4}
A.~Borodin and G.~Olshanski, \emph{Z--Measures on partitions,
Robinson--Schensted--Knuth correspondence, and $\beta=2$ random matrix
ensembles,} in: Random matrix models and their applications (P.~M.~Bleher and
A.~R.~Its, eds). Mathematical Sciences Research Institute Publications
\textbf{40}, Cambridge Univ. Press, 2001, 71--94; {\tt math/9905189}.



\bibitem{ol-BO5}
A.~Borodin and G.~Olshanski, \emph{Infinite random matrices and ergodic
measures,} Comm. Math. Phys., \textbf{223} (2001), 87--123; {\tt
math-ph/0010015}.



\bibitem{ol-BO6}
A.~Borodin and G.~Olshanski, \emph{Harmonic analysis on the
infinite--dimensional unitary group and determinantal point processes,} Ann.
Math, to appear; {\tt math/0109194}.


\bibitem{ol-BO7} A.~Borodin and G.~Olshanski, \emph{Z-measures on partitions
and their scaling limits,} 2002, {\tt math-ph/0210148}.



\bibitem{ol-DVJ}
D.~J.~Daley and D.~Vere-Jones, \emph{An introduction to the theory of point
processes,} Springer Series in Statistics, Springer, 1988.



\bibitem{ol-De}
P.~Deift, \emph{Integrable operators,} in: Differential operators and spectral
theory: M. Sh. Birman's 70th anniversary collection (V.~Buslaev, M.~Solomyak,
D.~Yafaev, eds.), American Mathematical Society Translations, ser. 2, vol. 189,
Providence, R.I., Amer. Math. Soc., 1999.


\bibitem{ol-DF} P. Diaconis and D. Freedman, \emph{Partial exchangeability and
sufficiency,} in: Statistics: Applications and New Directions (Calcutta, 1981),
Indian Statist. Inst., Calcutta, 1984, 205--236.


\bibitem{ol-Dy}
F.~J.~Dyson, \emph{Statistical theory of the energy levels of complex systems
I, II, III,} J. Math. Phys., \textbf{3} (1962), 140--156, 157--165, 166--175.




\bibitem{ol-Ed}
A. Edrei, \emph{On the generating functions of totally positive sequences II,}
J. Analyse Math., \textbf{2} (1952), 104--109.



\bibitem{ol-Ew1}
W.~J.~Ewens, \emph{The sampling theory of selectively neutral alleles,}
Theoret. Population Biology, \textbf{3}  (1972),   87--112.



\bibitem{ol-Ew2}
W.~J.~Ewens, \emph{Population Genetics Theory -- the Past and the Future,} in:
Mathematical and Statistical Developments of Evolutionary Theory (S.~Lessard,
ed.).  Proc. NATO ASI Symp., Kluwer, Dordrecht, 1990, 117--228.


\bibitem{ol-Gri} R.~C.~Griffiths, \emph{On the distribution of points in a Poisson Dirichlet
process,} J.~Appl. Probab. \textbf{25} (1988), 336--345.



\bibitem{ol-HN}
T.~Hida and H.~Nomoto, \emph{Gaussian measure on the projective limit space of
spheres,}  Proc. Japan Academy, \textbf{40} (1964), 31--34.


\bibitem{ol-Ig}
Ts. Ignatov, \emph{On a constant arising in the asymptotic theory of symmetric
groups and on Poisson--Dirichlet measures,}   Teor. Veroyatnost. Primenen.,
\textbf{27} (1982), no.~1, 129--140 (Russian); English translation: Theory
Probab. Appl. \textbf{27} (1982), 136--147.



\bibitem{ol-IIKS}
A.~R.~Its, A.~G.~Izergin, V.~E.~Korepin, N.~A.~Slavnov, \emph{Differential
equations for quantum correlation functions,}   Intern. J. Mod. Phys.,
\textbf{B4} (1990), 1003--1037.



\bibitem{ol-IO}
V.~Ivanov and G.~Olshanski, \emph{Kerov's central limit theorem for the
Plancherel measure on Young diagrams,}  in: Symmetric Functions 2001: Surveys
of Developments and Perspectives (S.~Fomin, ed.). NATO Science Series II.
Mathematics, Physics and Chemistry, vol. 74, Kluwer, 2001, 93--151.



\bibitem{ol-Jo}
K.~Johansson, \emph{Discrete orthogonal polynomial ensembles and the Plancherel
measure,}   Ann. Math. (2), \textbf{153} (2001),  no.~1, 259--296; {\tt
math/9906120}.



\bibitem{ol-Ke1}
S.~V.~Kerov, \emph{Combinatorial examples in the theory of AF-algebras,} in:
Differential geometry, Lie groups and mechanics X, Zapiski Nauchnykh Seminarov
LOMI \textbf{172} (1989), 55--67 (Russian); English translation: J.~Soviet
Math. \textbf{59} (1992), no.~5, 1063--1071.



\bibitem{ol-Ke2}
S.~V.~Kerov, \emph{Generalized Hall--Littlewood symmetric functions and
orthogonal polynomials,} in: Representation Theory and Dynamical Systems
(A.~M.~Vershik, ed.), Advances in Soviet Math., Vol. 9, Amer.  Math. Soc.,
Providence, R.I.,  1992, 67--94.



\bibitem{ol-Ke3}
S.~V.~Kerov, \emph{Gaussian limit for the Plancherel measure of the symmetric
group,}  Comptes Rendus Acad. Sci. Paris, S\'erie I, \textbf{316} (1993),
303--308.




\bibitem{ol-Ke4}
S.~V.~Kerov, \emph{Subordinators and the actions of permutations with
quasi--invariant measure,} in: Zapiski Nauchnyh Seminarov POMI \textbf{223}
(1995), 181--218 (Russian); English translation: J. Math. Sci.(New York)
\textbf{87} (1997), no.~6, 4094--4117.



\bibitem{ol-Ke5}
S.~V.~Kerov, \emph{Anisotropic Young diagrams and Jack symmetric functions,}
Funktsional. Anal. Prilozhen. \textbf{34} (2000), no.~1, 51--64 (Russian);
English translation: Funct. Anal. Appl., \textbf{34} (2000), no.~1, 41--51.



\bibitem{ol-KOO}
S.~Kerov, A.~Okounkov and G.~Olshanski, \emph{The boundary of Young graph with
Jack edge multiplicities,} Intern. Math. Res. Notices, \textbf{1998:4} (1998),
173--199; {\tt q-alg/9703037}.



\bibitem{ol-KO}
S.~Kerov and G.~Olshanski, \emph{Polynomial functions on the set of Young
diagrams,} Comptes Rendus Acad.\ Sci.\ Paris S\'er. I, \textbf{319} (1994),
121--126.


\bibitem{ol-KOV}
S.~Kerov, G.~Olshanski and A.~Vershik, \emph{Harmonic analysis on the infinite
symmetric group. A deformation of the regular representation,} Comptes Rendus
Acad. Sci. Paris, S\'er. I, \textbf{316} (1993), 773--778; detailed version in
preparation.



\bibitem{ol-KT}
S.~V.~Kerov and N.~V.~Tsilevich, \emph{Stick breaking process generates virtual
permutations with Ewens distribution,}  in: Zapiski Nauchnyh Seminarov POMI
\textbf{223} (1995), 162--180 (Russian); English translation: J. Math. Sci.
(New York), \textbf{87} (1997), no. 6, 4082--4093.



\bibitem{ol-Ki1}
J.~F.~C.~Kingman, \emph{Random discrete distributions,}  J.~Royal Statist. Soc.
B, \textbf{37} (1975), 1--22.



\bibitem{ol-Ki2}
    J.~F.~C.~Kingman,
\emph{The population structure associated with the Ewens sampling formula,}
Theoret. Population Biology, \textbf{11} (1977), 274--283.



\bibitem{ol-Ki3}
J.~F.~C.~Kingman, \emph{Random partitions in population genetics,
  Proc. Roy. Soc. London A.,} \textbf{361} (1978), 1--20.



\bibitem{ol-Ki4}
J.~F.~C.~Kingman, \emph{The representation of partition structures,}
  J.~London Math. Soc. (2), \textbf{18} (1978), 374--380.



\bibitem{ol-Ki5}
J.~F.~C.~Kingman, \emph{Poisson processes,} Oxford University Press, 1993.



\bibitem{ol-Len}
A.~Lenard, \emph{Correlation functions and the uniqueness of the state in
classical statistical mechanics,} Comm. Math. Phys, \textbf{30} (1973), 35--44.



\bibitem{ol-Ma1}
O.~Macchi, \emph{The coincidence approach to stochastic point processes,} Adv.
Appl. Prob., \textbf{7} (1975),   83--122.



\bibitem{ol-Ma2}
O.~Macchi, \emph{The fermion process --- a model of stochastic point process
with repulsive points,} in: Transactions of the Seventh Prague Conference on
Information Theory, Statistical Decision Functions, Random Processes and of the
Eighth European Meeting of Statisticians (Tech. Univ. Prague, Prague, 1974),
Vol. A, Reidel, Dordrecht, 1977,   391--398.


\bibitem{ol-MvN} F.~J.~Murray and J.~von~Neumann, \emph{On rings of operators
IV,}  Ann. Math. \textbf{44} (1943), 716--808.




\bibitem{ol-Me}
M.~L.~Mehta, \emph{Random matrices,}  2nd edition, Academic Press, 1991.


\bibitem{ol-Ne}
Yu. A. Neretin, \emph{Hua type integrals over unitary groups and over
projective limits of unitary groups,} Duke Math. J., \textbf{114} (2002),
239--266; {\tt math-ph/0010014}.




\bibitem{ol-NW}
T.~Nagao, M.~Wadati, \emph{Correlation functions of random matrix ensembles
related to classical orthogonal polynomials,}  J. Phys. Soc. Japan, \textbf{60}
(1991), no.~10, 3298--3322.



\bibitem{ol-Na}
M.~A.~Naimark, \emph{Normed rings,}  Nauka, Moscow, 1962 (Russian); English
translation: \emph{Normed algebras,}  Wolters--Noordhoff, Groningen, 1972.



\bibitem{ol-Ok1}
A.~Yu.~Okounkov, \emph{Thoma's theorem and representations of infinite
bisymmetric group,} Funktsion. Anal. Prilozhen. \textbf{28} (1994), no.~2,
31--40 (Russian); English translation: Funct. Anal. Appl., \textbf{28} (1994),
no. 2, 101--107.

\bibitem{ol-Ok1a}
A.~Yu.~Okounkov \emph{On representations of the infinite symmetric group,} in:
Representation Theory, Dynamical Systems, Combinatorial and Algorithmic Methods
II,  Zap.\ Nauchn.\ Semin.\ POMI  (A.~M.~Vershik, ed.) \textbf{240} (1997),
167--229 (Russian); English translation: J.~Math. Sci. (New York), \textbf{96}
(1999), no.~5, 3550--3589.



\bibitem{ol-Ok2}
A.~Okounkov, \emph{$SL(2)$ and z--measures,} in: Random matrix models and their
applications (P.~M.~Bleher and A.~R.~Its, eds). Mathematical Sciences Research
Institute Publications \textbf{40}, Cambridge Univ. Press, 2001, 407--420; {\tt
math/0002136}.



\bibitem{ol-Ok3}
A.~Okounkov, \emph{Infinite wedge and measures on partitions,}  Selecta Math.
(New Ser.) \textbf{7} (2001), 57--81; {\tt math/9907127}.



\bibitem{ol-OO}
A.~Okounkov and G.~Olshanski, \emph{Shifted Schur functions,}  Algebra i
Analiz, \textbf{9} (1997), no.~2, 73--146 (Russian); English translation:
St.~Petersburg Math. J. \textbf{9} (1998), no.~2, 239--300.



\bibitem{ol-Ol1}
G.~Olshanski, \emph{Unitary representations of infinite-dimensional pairs
$(G,K)$ and the formalism of R.~Howe,} Doklady AN SSSR, \textbf{269} (1983),
33--36 (Russian); English translation: Soviet Math. Doklady, \textbf{27}
(1983), no.~2, 290--294.



\bibitem{ol-Ol2}
G.\ Olshanski, \emph{Unitary representations of infinite--dimensional pairs
$(G,K)$ and the formalism of R.~Howe,} in: Representation of Lie Groups and
Related Topics (A.~Vershik and D.~Zhelobenko, eds.), Advanced Studies in
Contemporary Math. \textbf{7}, Gordon and Breach Science Publishers, New York
etc., 1990, 269--463.



\bibitem{ol-Ol3}
G.\ Olshanski, \emph{Unitary representations of $(G,K)$-pairs connected with
the infinite symmetric group $S(\infty)$,}  Algebra i Analiz, \textbf{1}
(1989), no. 4, 178--209 (Russian); English translation: Leningrad Math. J.
\textbf{1} (1990), 983--1014.



\bibitem{ol-Ol4}
G.\ Olshanski, \emph{The problem of harmonic analysis on the
infinite--dimensional unitary group,}  J. Funct. Anal., to appear; {\tt
math/0109193}.



\bibitem{ol-ORV}
G.~Olshanski, A.~Regev and A.~Vershik, \emph{Frobenius--Schur functions,} in:
Studies in Memory of Issai Schur (A.~Joseph, A.~Melnikov, R.~Rentschler, eds.),
Birkh\"auser, to appear; {\tt math/0110077}.




\bibitem{ol-P.I}
G.~Olshanski, \emph{Point processes and the infinite symmetric group. Part I:
The general formalism and the density function,} 1998, {\tt math/9804086}.



\bibitem{ol-P.II}
A.~Borodin, \emph{Point processes and the infinite symmetric group. Part II:
Higher correlation functions,} 1998, {\tt math/9804087}.



\bibitem{ol-P.III}
A.~Borodin and G.~Olshanski, \emph{Point processes and the infinite symmetric
group. Part III: Fermion point processes,} 1998, {\tt math/9804088}.



\bibitem{ol-P.IV}
A.~Borodin, \emph{Point processes and the infinite symmetric group. Part IV:
Matrix Whittaker kernel,} 1998, {\tt math/9810013}.



\bibitem{ol-P.V} G.~Olshanski, \emph{Point processes and the infinite symmetric
group. Part V: Analysis of the matrix Whittaker kernel,} 1998, {\tt
math/9810014}.



\bibitem{ol-Pi}
D.~Pickrell, \emph{Measures on infinite dimensional Grassmann manifold,}
J.~Funct. Anal., \textbf{70} (1987), 323--356.



\bibitem{ol-R}
N.~A.~Rozhkovskaya, \emph{Multiplicative distributions on Young graph},  in:
Representation Theory, Dynamical Systems, Combinatorial and Algorithmic Methods
II (A.~M.~Vershik, ed.), Zapiski Nauchnykh Seminarov POMI {\textbf{240}, Nauka,
St.~Petersburg,   1997,   246-257 (Russian); English translation: J. Math. Sci.
(New York) 96 (1999), no.~5, 3600--3608



\bibitem{ol-Shi}
H.~Shimomura, \emph{On the construction of invariant measure over the
orthogonal group on the Hilbert space by the method of Cayley transformation,}
Publ. RIMS Kyoto Univ., \textbf{10}  (1974/75),   413--424.



\bibitem{ol-So}
A.~Soshnikov, \emph{Determinantal random point fields,} Uspekhi Mat. Nauk
\textbf{55} (2000), no.~5, 107--160 (Russian); English translation: Russian
Math. Surveys, \textbf{55} (2000), no.~5, 923--975; {\tt math/0002099}.



\bibitem{ol-SV}
S.~Stratila and D.~Voiculescu, \emph{Representations of AF--algebras and of the
group $U(\infty)$,} Springer Lecture Notes \textbf{486}, 1975.



\bibitem{ol-Th1}
E.~Thoma, \emph{Die unzerlegbaren, positive--definiten Klassenfunktionen der
abz\"ahlbar unendlichen, symmetrischen Gruppe,} Math. Zeitschr., \textbf{85}
(1964),   40--61.



\bibitem{ol-Th2}
E.~Thoma, \emph{Characters of infinite groups,} in: Operator Algebras and Group
Representations (Gr. Arsene, S. Str{\u a}til{\u a}, A. Verona, and D.
Voiculescu, eds), Vol.~2, Pitman,   1984,  pp. 23--32.



\bibitem{ol-TW}
C.~A.~Tracy and H.~Widom, \emph{Universality of distribution functions of
random matrix theory II,}   in: Integrable Systems: From Classical to Quantum
(J.~Harnad, G.~Sabidussi, and P.~Winternitz, eds). CRM Proceedings \& Lecture
Notes \textbf{26}. Amer. Math. Soc., Providence, 2000, pp. 251--264.



\bibitem{ol-Ve1}
A.~M.~Vershik, \emph{Description of invariant measures for the actions of some
infinite-dimensional groups,} Doklady AN SSSR, \textbf{218} (1974), 749--752
(Russian); English translation: Soviet Math. Doklady, \textbf{15} (1974),
1396--1400.



\bibitem{ol-Ve2}
A.~M.~Vershik, \emph{Asymptotic distribution of decompositions of natural
numbers into prime divisors,}  Dokl. Akad. Nauk SSSR, \textbf{289}  (1986),
no.~2,   269--272; English translation: Soviet Math. Doklady, \textbf{34}
(1987), 57--61.



\bibitem{ol-Ve3}
A.~M.~Vershik, \emph{Statistical mechanics of combinatorial partitions, and
their limit shapes,}  Funktsional. Anal.  Prilozhen. \textbf{30} (1996), no.~2,
19--39 (Russian); English translation: Funct. Anal. Appl., \textbf{30}  (1996),
90--105.



\bibitem{ol-VK1}
A.~M.~Vershik and S.~V.~Kerov, \emph{Characters and factor representations of
the infinite symmetric group,}   Doklady AN SSSR, \textbf{257} (1981),
1037--1040 (Russian); English translation: Soviet Math. Doklady \textbf{23}
(1981), 389--392.



\bibitem{ol-VK2}
A.~M.~Vershik and S.~V.~Kerov, \emph{Asymptotic theory of characters of the
symmetric group,} Funktsion. Anal. Prilozhen.,  \textbf{15} (1981), no.~4,
15--27 (Russian); English translation: Funct. Anal. Appl., \textbf{15} (1981),
no.~4, 246--255.



\bibitem{ol-VK3}
A.~M.~Vershik and S.~V.~Kerov, \emph{Characters and factor representations of
the infinite unitary group,} Doklady AN SSSR, \textbf{267} (1982), no.~2,
272--276 (Russian); English translation: Soviet Math. Doklady, \textbf{26}
(1982), 570--574.



\bibitem{ol-VK4}
A.~M.~Vershik and S.~V.~Kerov, \emph{Locally semisimple algebras. Combinatorial
theory and the $K_0$ functor,} in: Itogi Nauki, Sovr. Probl. Mat., Noveish.
Dostizh., VINITI, \textbf{26} (1985), 3--56 (Russian); English translation:
J.~Soviet Math., \textbf{38} (1987), 1701--1733.



\bibitem{ol-VS}
A.~M.~Vershik and A.~A.~Shmidt, \emph{Limit measures arising in the asymptotic
theory of symmetric groups I, II},   Teor. Veroyatnost.  Primenen., \textbf{22}
(1977),   no.~1, 72--88, \textbf{23}  (1978),   no.~1, 42--54 (Russian);
English translation: Theory of Probab.  Appl. \textbf{22} (1977), 70-85,
\textbf{23} (1978), 36--49.




\bibitem{ol-Vo}
D.~Voiculescu, \emph{Repr\'esentations factorielles de type {\rm II}${}_1$ de
$U(\infty)$,}   J.\ Math.\ Pures  Appl., \textbf{55}  (1976),   1--20 .



\bibitem{ol-Was}
A.~J.~Wassermann, \emph{Automorphic actions of compact groups on operator
algebras,}  Thesis, University of Pennsylvania, 1981.


\bibitem{ol-Wat}
G.~A.~Watterson, \emph{The stationary distribution of the infinitely
many--alleles diffusion model,}   J.~Appl. Probab., \textbf{13}  (1976),
639--651.



\bibitem{ol-Ya1}
Y.~Yamasaki, \emph{Projective limit of Haar measures on $O(n)$,} Publ. RIMS,
Kyoto Univ.,} \textbf{8}  (1972/73),   141--149.



\bibitem{ol-Ya2}
Y.~Yamasaki, \emph{Kolmogorov's extension theorem for infinite measure,} Publ.
RIMS, Kyoto Univ., \textbf{10}  (1974/75),   381--411.



\end{thebibliography}
\end{document}